\documentclass[a4paper]{article}

\usepackage{amsfonts,amssymb,amsthm}

\makeatletter
\newif\iflabel\labelfalse \let\w@label=\label
\def\label{\global\labeltrue\w@label}
\let\w@eqn=\equation
\def\equation{\global\labelfalse\w@eqn}
\def\endequation{\iflabel\@eeq\else\@eeqw\fi$$\global\@ignoretrue}
\def\@eeq{\eqno \@eqnnum}
\def\@eeqw{\addtocounter{equation}{-1}}
\newskip\bw@\newskip\bws@
\def\dc@@pq{\global \bw@=\belowdisplayskip \global\bws@=\belowdisplayshortskip
   \global\belowdisplayshortskip=3pt plus -3pt
   \global\belowdisplayskip=\belowdisplayshortskip
   \hskip 500pt minus 500pt\relax}
\def\dc@pq{\dc@@pq$$}
\def\fc@pq{$$\hskip 500pt minus 500pt\global
   \belowdisplayshortskip=\bws@\global\belowdisplayskip=\bw@}
\def\coupeq{\dc@pq\fc@pq}
\def\coupepas{\par\noindent\begin{minipage}{\textwidth}}
\def\jusquela{\end{minipage}}
\makeatother

\catcode`\@=11\relax
\newwrite\@unused
\def\typeout#1{{\let\protect\string\immediate\write\@unused{#1}}}
\def\@nnil{\@nil}
\def\@empty{}
\def\@psdonoop#1\@@#2#3{}
\def\@psdo#1:=#2\do#3{\edef\@psdotmp{#2}\ifx\@psdotmp\@empty \else
    \expandafter\@psdoloop#2,\@nil,\@nil\@@#1{#3}\fi}
\def\@psdoloop#1,#2,#3\@@#4#5{\def#4{#1}\ifx #4\@nnil \else
       #5\def#4{#2}\ifx #4\@nnil \else#5\@ipsdoloop #3\@@#4{#5}\fi\fi}
\def\@ipsdoloop#1,#2\@@#3#4{\def#3{#1}\ifx #3\@nnil 
       \let\@nextwhile=\@psdonoop \else
      #4\relax\let\@nextwhile=\@ipsdoloop\fi\@nextwhile#2\@@#3{#4}}
\def\@tpsdo#1:=#2\do#3{\xdef\@psdotmp{#2}\ifx\@psdotmp\@empty \else
    \@tpsdoloop#2\@nil\@nil\@@#1{#3}\fi}
\def\@tpsdoloop#1#2\@@#3#4{\def#3{#1}\ifx #3\@nnil 
       \let\@nextwhile=\@psdonoop \else
      #4\relax\let\@nextwhile=\@tpsdoloop\fi\@nextwhile#2\@@#3{#4}}
\def\psdraft{
	\def\@psdraft{0}
}
\def\psfull{
	\def\@psdraft{100}
}
\psfull
\newif\if@prologfile
\newif\if@postlogfile
\newif\if@noisy
\def\pssilent{
	\@noisyfalse
}
\def\psnoisy{
	\@noisytrue
}
\psnoisy
\newif\if@bbllx
\newif\if@bblly
\newif\if@bburx
\newif\if@bbury
\newif\if@height
\newif\if@width
\newif\if@rheight
\newif\if@rwidth
\newif\if@clip
\newif\if@verbose
\def\@p@@sclip#1{\@cliptrue}
\def\@p@@sfile#1{
		   \def\@p@sfile{#1}
}
\def\@p@@sfigure#1{\def\@p@sfile{#1}}
\def\@p@@sbbllx#1{
		\@bbllxtrue
		\dimen100=#1
		\edef\@p@sbbllx{\number\dimen100}
}
\def\@p@@sbblly#1{
		\@bbllytrue
		\dimen100=#1
		\edef\@p@sbblly{\number\dimen100}
}
\def\@p@@sbburx#1{
		\@bburxtrue
		\dimen100=#1
		\edef\@p@sbburx{\number\dimen100}
}
\def\@p@@sbbury#1{
		\@bburytrue
		\dimen100=#1
		\edef\@p@sbbury{\number\dimen100}
}
\def\@p@@sheight#1{
		\@heighttrue
		\dimen100=#1
   		\edef\@p@sheight{\number\dimen100}
}
\def\@p@@swidth#1{
		\@widthtrue
		\dimen100=#1
		\edef\@p@swidth{\number\dimen100}
}
\def\@p@@srheight#1{
		\@rheighttrue
		\dimen100=#1
		\edef\@p@srheight{\number\dimen100}
}
\def\@p@@srwidth#1{
		\@rwidthtrue
		\dimen100=#1
		\edef\@p@srwidth{\number\dimen100}
}
\def\@p@@ssilent#1{ 
		\@verbosefalse
}
\def\@p@@sprolog#1{\@prologfiletrue\def\@prologfileval{#1}}
\def\@p@@spostlog#1{\@postlogfiletrue\def\@postlogfileval{#1}}
\def\@cs@name#1{\csname #1\endcsname}
\def\@setparms#1=#2,{\@cs@name{@p@@s#1}{#2}}
\def\ps@init@parms{
		\@bbllxfalse \@bbllyfalse
		\@bburxfalse \@bburyfalse
		\@heightfalse \@widthfalse
		\@rheightfalse \@rwidthfalse
		\def\@p@sbbllx{}\def\@p@sbblly{}
		\def\@p@sbburx{}\def\@p@sbbury{}
		\def\@p@sheight{}\def\@p@swidth{}
		\def\@p@srheight{}\def\@p@srwidth{}
		\def\@p@sfile{}
		\def\@p@scost{10}
		\def\@sc{}
		\@prologfilefalse
		\@postlogfilefalse
		\@clipfalse
		\if@noisy
			\@verbosetrue
		\else
			\@verbosefalse
		\fi
}
\def\parse@ps@parms#1{
	 	\@psdo\@psfiga:=#1\do
		   {\expandafter\@setparms\@psfiga,}}
\newif\ifno@bb
\newif\ifnot@eof
\newread\ps@stream
\def\bb@missing{
	\if@verbose{
		\typeout{psfig: searching \@p@sfile \space  for bounding box}
	}\fi
	\openin\ps@stream=\@p@sfile
	\no@bbtrue
	\not@eoftrue
	\catcode`\%=12
	\loop
		\read\ps@stream to \line@in
		\global\toks200=\expandafter{\line@in}
		\ifeof\ps@stream \not@eoffalse \fi
		\@bbtest{\toks200}
		\if@bbmatch\not@eoffalse\expandafter\bb@cull\the\toks200\fi
	\ifnot@eof \repeat
	\catcode`\%=14
}	
\catcode`\%=12
\newif\if@bbmatch
\def\@bbtest#1{\expandafter\@a@\the#1
\long\def\@a@#1
\long\def\bb@cull#1 #2 #3 #4 #5 {
	\dimen100=#2 bp\edef\@p@sbbllx{\number\dimen100}
	\dimen100=#3 bp\edef\@p@sbblly{\number\dimen100}
	\dimen100=#4 bp\edef\@p@sbburx{\number\dimen100}
	\dimen100=#5 bp\edef\@p@sbbury{\number\dimen100}
	\no@bbfalse
}
\catcode`\%=14
\def\compute@bb{
		\no@bbfalse
		\if@bbllx \else \no@bbtrue \fi
		\if@bblly \else \no@bbtrue \fi
		\if@bburx \else \no@bbtrue \fi
		\if@bbury \else \no@bbtrue \fi
		\ifno@bb \bb@missing \fi
		\ifno@bb \typeout{FATAL ERROR: no bb supplied or found}
			\no-bb-error
		\fi
		\count203=\@p@sbburx
		\count204=\@p@sbbury
		\advance\count203 by -\@p@sbbllx
		\advance\count204 by -\@p@sbblly
		\edef\@bbw{\number\count203}
		\edef\@bbh{\number\count204}
}
\def\in@hundreds#1#2#3{\count240=#2 \count241=#3
		     \count100=\count240	
		     \divide\count100 by \count241
		     \count101=\count100
		     \multiply\count101 by \count241
		     \advance\count240 by -\count101
		     \multiply\count240 by 10
		     \count101=\count240	
		     \divide\count101 by \count241
		     \count102=\count101
		     \multiply\count102 by \count241
		     \advance\count240 by -\count102
		     \multiply\count240 by 10
		     \count102=\count240	
		     \divide\count102 by \count241
		     \count200=#1\count205=0
		     \count201=\count200
			\multiply\count201 by \count100
		 	\advance\count205 by \count201
		     \count201=\count200
			\divide\count201 by 10
			\multiply\count201 by \count101
			\advance\count205 by \count201
		     \count201=\count200
			\divide\count201 by 100
			\multiply\count201 by \count102
			\advance\count205 by \count201
		     \edef\@result{\number\count205}
}
\def\compute@wfromh{
		\in@hundreds{\@p@sheight}{\@bbw}{\@bbh}
		\edef\@p@swidth{\@result}
}
\def\compute@hfromw{
		\in@hundreds{\@p@swidth}{\@bbh}{\@bbw}
		\edef\@p@sheight{\@result}
}
\def\compute@handw{
		\if@height 
			\if@width
			\else
				\compute@wfromh
			\fi
		\else 
			\if@width
				\compute@hfromw
			\else
				\edef\@p@sheight{\@bbh}
				\edef\@p@swidth{\@bbw}
			\fi
		\fi
}
\def\compute@resv{
		\if@rheight \else \edef\@p@srheight{\@p@sheight} \fi
		\if@rwidth \else \edef\@p@srwidth{\@p@swidth} \fi
}
\def\compute@sizes{
	\compute@bb
	\compute@handw
	\compute@resv
}
\def\psfig#1{\vbox {
	\ps@init@parms
	\parse@ps@parms{#1}
	\compute@sizes
	\ifnum\@p@scost<\@psdraft{
		\if@verbose{
			\typeout{psfig: including \@p@sfile \space }
		}\fi
		\special{ps::[begin] 	\@p@swidth \space \@p@sheight \space
				\@p@sbbllx \space \@p@sbblly \space
				\@p@sbburx \space \@p@sbbury \space
				startTexFig \space }
		\if@clip{
			\if@verbose{
				\typeout{(clip)}
			}\fi
			\special{ps:: doclip \space }
		}\fi
		\if@prologfile
		    \special{ps: plotfile \@prologfileval \space } \fi
		\special{ps: plotfile \@p@sfile \space }
		\if@postlogfile
		    \special{ps: plotfile \@postlogfileval \space } \fi
		\special{ps::[end] endTexFig \space }
		\vbox to \@p@srheight true sp{
			\hbox to \@p@srwidth true sp{
				\hss
			}
		\vss
		}
	}\else{
		\vbox to \@p@srheight true sp{
		\vss
			\hbox to \@p@srwidth true sp{
				\hss
				\if@verbose{
					\@p@sfile
				}\fi
				\hss
			}
		\vss
		}
	}\fi
}}
\catcode`\@=12\relax

\newcommand{\vs}{\vspace{0.3cm}}
\newcommand{\dr}{\partial}
\newcommand{\R}{{\bf R}}
\newcommand{\N}{{\bf N}}
\newcommand{\Z}{{\bf Z}}
\newcommand{\II}{I\hspace{-0.1cm}I}
\newcommand{\III}{I\hspace{-0.1cm}I\hspace{-0.1cm}I}
\newcommand{\tr}{\mbox{\rm tr}}
\newcommand{\ric}{\mbox{\rm ric}}
\newcommand{\cotg}{\mbox{\rm cotg}}
\newcommand{\SO}{\mbox{\rm SO}}
\newcommand{\ricb}{\overline{\ric}}
\newcommand{\Rc}{\check{R}}
\newcommand{\diam}{\mbox{\rm diam}}
\newcommand{\area}{\mbox{\rm area}}

\def\pointir{\unskip  {. --- \ignorespaces }\hskip0cm}

\newcommand{\deltab}{\overline{\delta}}

\newcommand{\gba}{\overline{g}}
\newcommand{\hb}{\overline{h}}
\newcommand{\Db}{\overline{D}}
\newcommand{\Rb}{\overline{R}}
\newcommand{\Kb}{\overline{K}}
\newcommand{\Sb}{\overline{S}}

\newcommand{\pt}{\tilde{p}}
\newcommand{\ut}{\tilde{u}}
\newcommand{\yt}{\tilde{y}}
\newcommand{\Bt}{\tilde{B}}
\newcommand{\Ft}{\tilde{F}}
\newcommand{\Kt}{\tilde{K}}
\newcommand{\ept}{\tilde{\epsilon}}
\newcommand{\gat}{\tilde{\gamma}}
\newcommand{\kappat}{\tilde{\kappa}}
\newcommand{\thetat}{\tilde{\theta}}

\newcommand{\rhot}{\tilde{\rho}}
\newcommand{\taut}{\tilde{\tau}}
\newcommand{\phit}{\tilde{\phi}}
\newcommand{\sit}{\tilde{\sigma}}
\newcommand{\Phit}{\tilde{\Phi}}
\newcommand{\Sigmat}{\tilde{\Sigma}}
\newcommand{\Omt}{\tilde{\Omega}}

\newcommand{\kad}{\stackrel{\bullet}{\kappa}}
\newcommand{\gd}{\stackrel{\bullet}{g}}

\newtheorem{prop}{Proposition}[section]
\newtheorem{df}[prop]{Definition}
\newtheorem{lemma}[prop]{Lemma}
\newtheorem{thm}[prop]{Theorem}
\newtheorem{cor}[prop]{Corollary}
\newtheorem{asser}[prop]{Assertion}
\newtheorem{remark}[prop]{Remark}

\newenvironment{thn}[1]{\vskip 0.2cm \noindent{\bf Theorem #1.} \it}{\rm
\vspace{0.2cm}} 
\newenvironment{lmn}[1]{\vskip 0.2cm \noindent{\bf Lemma #1.} \it}{\rm
\vspace{0.2cm}} 

\newcommand{\btm}{\begin{thm}}
\newcommand{\etm}{\end{thm}}
\newcommand{\blm}{\begin{lemma}}
\newcommand{\elm}{\end{lemma}}
\newcommand{\bcr}{\begin{cor}}
\newcommand{\ecr}{\end{cor}}
\newcommand{\bdf}{\begin{df}}
\newcommand{\edf}{\end{df}}
\newcommand{\bprop}{\begin{prop}}
\newcommand{\eprop}{\end{prop}}
\newcommand{\bas}{\begin{asser}}
\newcommand{\eas}{\end{asser}}
\newcommand{\beq}{\begin{equation}}
\newcommand{\eeq}{\end{equation}}
\newcommand{\bpv}{\begin{proof}}
\newcommand{\epv}{\end{proof}}
\newcommand{\bit}{\begin{itemize}}
\newcommand{\eit}{\end{itemize}}
\newcommand{\bpn}{\begin{pfn}}
\newcommand{\epn}{\end{pfn}}
\newcommand{\btn}{\begin{thn}}
\newcommand{\etn}{\end{thn}}
\newcommand{\bln}{\begin{lmn}}
\newcommand{\eln}{\end{lmn}}

\newenvironment{pfn}[1]{\vskip 0.2cm \noindent{\it Proof #1.}}{$\square$
\vspace{0.2cm}} 

\newcommand{\Omb}{\overline{\Omega}}
\newcommand{\Sib}{\overline{\Sigma}}
\newcommand{\gb}{\overline{g}}
\newcommand{\Ub}{\overline{U}}
\newcommand{\Wb}{\overline{W}}
\newcommand{\db}{\overline{\partial}}

\newcommand{\Met}{\mathcal{M}et}
\newcommand{\Imm}{\mathcal{I}mm}
\newcommand{\CMet}{\mathcal{CM}et}
\newcommand{\cE}{\mathcal{E}}
\newcommand{\cM}{\mathcal{M}}
\newcommand{\cS}{\mathcal{S}}
\newcommand{\cP}{\mathcal{P}}
\newcommand{\CImm}{\mathcal{CI}mm}
\newcommand{\gab}{\overline{\gamma}}
\newcommand{\hyp}{\mathbf{H}^3}
\newcommand{\dhyp}{\partial\hyp}
\newcommand{\cL}{\mathcal{L}}
\newcommand{\isom}{\mathrm{Isom}}
\newcommand{\im}{\mathrm{Im}}
\newcommand{\Na}{\nabla}
\newcommand{\Nat}{\tilde{\nabla}}
\newcommand{\Sit}{\tilde{\Sigma}}

\newcommand{\eps}{\epsilon}
\newcommand{\ga}{\gamma}
\newcommand{\si}{\sigma}
\newcommand{\om}{\omega}
\newcommand{\Ga}{\Gamma}
\newcommand{\La}{\Lambda}
\newcommand{\Si}{\Sigma}
\newcommand{\Om}{\Omega}

\begin{document}

\title{Complete surfaces with negative extrinsic curvature}

\date{December 1999}

\author{Jean-Marc Schlenker\thanks{
Math{\'e}matiques, UMR 8628 du CNRS, B{\^a}t. 425, Uni\-ver\-sit{\'e}
Paris-Sud, F-91405 Orsay Cedex, France~; currently: FIM, ETHZ,
R{\"a}mistr. 101, CH-8092 Z{\"u}rich. 
\texttt{jean-marc.schlenker@math.u-psud.fr}.
\texttt{http://www.math.u-psud.fr/\~{ }schlenker}.
} 
}

\maketitle

\begin{abstract}

N. V. Efimov \cite{Ef1} proved that there is no complete, smooth surface
in $\R^3$ with uniformly negative curvature. We extend this to isometric
immersions in a 3-manifold with pinched curvature: if $M^3$ has sectional
curvature between two constants $K_2$ and $K_3$, then there exists
$K_1 < \min(K_2, 0)$ such that $M$ contains no smooth, complete immersed
surface with curvature below $K_1$. Optimal values of $K_1$ are
determined. This results rests on a phenomenon of propagations for
degenerations of solutions of hyperbolic Monge-Amp{\`e}re equations.

\bigskip

\begin{center} {\bf R{\'e}sum{\'e}} \end{center}

N. V. Efimov \cite{Ef1} a montr{\'e} qu'il n'existe pas de surface compl{\`e}te
{\`a} courbure uniform{\'e}ment n{\'e}gative dans $\R^3$. On {\'e}tend ce r{\'e}sultat aux
immersions isom{\'e}triques dans les 3-vari{\'e}t{\'e}s {\`a} courbure pinc{\'e}e: si $M^3$
a sa courbure sectionnelle comprise entre deux constantes $K_2$ et
$K_3$, alors il existe une constante $K_1<\min(K_2, 0)$ telle que $M$ ne
contient pas de surface immerg{\'e}e compl{\`e}te et r{\'e}guli{\`e}re {\`a} courbure
inf{\'e}rieure {\`a} $K_1$. Des valeurs optimales de $K_1$ sont d{\'e}termin{\'e}es. Ce
r{\'e}sultat repose sur un ph{\'e}nom{\`e}ne de propagation pour les d{\'e}g{\'e}n{\'e}rescences
de solutions d'{\'e}quations de Monge-Amp{\`e}re hyperboliques.
 
\end{abstract}

\maketitle

\vspace{0.4cm}

{\bf AMS classifications:} 53C45, 58G16, 35L55.

\vspace{0.4cm}

{\bf Key-words:} isometric, immersion, surface, Monge-Amp{\`e}re,
hyperbolic.

\bigskip

Hilbert \cite{Hil} proved that there is no smooth isometric immersion of
the hyperbolic plane $H^2$ into the Euclidean 3-space $\R^3$. This was
extended by Efimov, who replaced $H^2$ by any complete surface with
uniformly negative curvature:

\begin{thm}[N. V. Efimov \cite{Ef1}] \label{efimov}
Let $(\Si, \si)$ be a smooth, complete Riemannian surface with curvature
$K\leq -1$. Then $(\Si, \si)$ has no $C^2$ isometric immersion into
$\R^3$. 
\end{thm}

This result was proved using some subtle geometric constructions,
strongly based on the Euclidean structure of the target space. More
details can be found in \cite{Ef6}, \cite{KM} or in \cite{BS,Roz3}, and
some extensions and related results in \cite{Ef2,Ef3,Ef5}.

It seems rather natural to try to extend Hilbert's result further by
replacing also $\R^3$ by a Riemannian manifold. This was started in
\cite{efj}, where the target space can be a Riemannian or Lorentzian
3-dimensional space-form. The present paper treats the case where it is
a Riemannian manifold with pinched curvature.

\begin{thm} \label{thm-a}
Let $(M,\mu)$ be a complete Riemannian 3-manifold, with sectional
curvature $K_M$ between two constants $K_2 \leq K_3$. Let 
$(\Si,\si)$ be a complete Riemannian surface, with curvature $K_{\Si}\leq
K_1$, with $K_1<0$, $K_1<K_2 \leq K_3$, and:
\bit
\item either $K_3\geq 0$ and
\beq (K_3-K_2)^2 < 16 |K_1| (K_2-K_1)~; \eeq
\item or $K_3\leq 0$ and
\beq (K_3-K_2)^2 < 16 (K_3-K_1)(K_2-K_1)~. \eeq
\eit
Finally, suppose that $\| \nabla (K_{\Si}^{-1/2})\|$ and $\| \nabla
K_{M}\|$ are bounded.
Then there exists no $C^{3}$ isometric immersion from $(\Si,\si)$ into
$(M,\mu)$.
\end{thm}

The meaning of ``$\| \nabla K_{M}\|$ bounded'' demands some
precisions. Let $m\in M$, let $P$ be a 2-plane in $T_mM$, and let
$c:[0,1]\rightarrow M$ be a smooth curve with $c(0)=m$. For $t\in [0,
1]$, call $P_t$ the parallel transport of $P$ at $c(t)$ along $c([0,
t])$, and let $K(t)$ be the sectional curvature of $M$ on $P_t$. Then our
hypothesis is that $| K'(t)|$ is bounded by some fixed constant.

\bigskip

The proof of theorem \ref{thm-a} rests on two ideas, one of a geometric
and the other of an analytical nature. 

The geometric point concerns which
objects, induced on a surface by an immersion, are to be considered. Of
course, one could consider the induced metric -- also called the first
fundamental form $I$ of the immersion -- along with its Levi-Civita
connection $\Na$ and the ``Weingarten operator'' $B$, which satisfies
what can be described as a Monge-Amp{\`e}re equation of hyperbolic type:
$\det(B)$ is equal to the extrinsic curvature of the immersion (which is
negative here), while $d^{\Na}B$ is equal to another term given by the
Coddazi equation, which is bounded. There are some ``dual'' objects,
however, which are of greater use: the third fundamental form
$\III$ of the surface, and the inverse $\Bt$ of $B$. The ``new'' point
is that the ``right'' connection to use is not the Levi-Civita connection
of $\III$, but rather another connection, called $\Nat$, which is
compatible with $\III$ and has bounded torsion. $\Bt$ then satisfies a
very simple equation: $\det(\Bt)$ is 
again given by the extrinsic curvature, while $d^{\Nat}\Bt=0$. When the
ambiant space 
has constant curvature, $\Nat$ is indeed the Levi-Civita connection of
$\III$. 

The analytical fact which is important in the proof is about
propagations of degenerations of sequences of solutions of some
hyperbolic Monge-Amp{\`e}re equations. Remember again that isometric
immersions of surfaces are described analytically as solutions 
of Monge-Amp{\`e}re equations. When the extrinsic curvature of the immersed
surface is positive, the equations are elliptic, and this case is rather
well understood \cite{Po,CNS,CNS3,CNS4,CNS5,L1,these,iie}.  
A fundamental point is
that solutions of those equations have no isolated singularities:
rather, if a sequence of solution has a limit which is degenerate at a
point, then (for some subsequence) the same happens along a
geodesic. This phenomenon has been 
studied completely by F. Labourie in \cite{L2,L1,L-MA} (see \cite{BK}
for some related problems). It is
interesting to remark that, for complex Monge-Amp{\`e}re solutions, the
geometric nature of the locus of degeneration of sequences of solutions
also plays a major role (see e.g. \cite{Nadel}).

On the other hand, it has been knows since \cite{Roz1} that surfaces
with negative curvature in $\R^3$ can have an isolated
singularity. Nonetheless, a
phenomenon of propagation of degenerations of sequences of solutions of
hyperbolic Monge-Amp{\`e}re equations appears when the singularities are
supposed to be bad enough. Here is an example of such a result.

\begin{thm}[\cite{efj}] \label{thpr}
Let $D$ be a disk with a smooth Riemannian metric $g$ with curvature
$K<-1$, and let $(\phi_n)_{n\in \N}$ be a sequence of isometric
immersions of
$(D,g)$ into $\R^3$. Let $x_0\in D$, and let $y_0\in \R^3$ be such that,
for all $n$, $\phi_n(x_0)=y_0$. Suppose that $(\phi_n)$ is degenerate at
$x_0$, in the sense that there exists a geodesic segment $\ga_0$ with
$\ga_0(0)=x_0$ such that:
$$ \forall \eps>0, \exists n\in \N,
\int_0^\eps (\III_n(\ga'_0(s),\ga'_0(s))^{1/2} ds \geq \frac{1}{\eps}~. $$
Then there exists a subsequence $(\psi_n)_{n\in \N}$ of
$(\phi_n)_{n\in \N}$ and a maximal geodesic segment $g$ going through
$x_0$ such that $(\psi_n)$ is degenerate along $g$~:
$$ \forall \eps>0, \exists N\in \N, \forall n\geq N, \forall x\in g,
\exists y\in B_{\mu}(x, \eps), H_n(y)\geq 1/\eps~. $$
Moreover $(\psi_{n|g})_{n\in \N}$ converges $C^0$ towards an isometry
from $g$ to a geodesic segment of $\R^3$.
\end{thm}

This kind of propagation is essentially responsible for a crucial point
of the proof of theorem \ref{thm-a}, namely that $(\Si, \III, \Nat)$ is
``convex'' in a precise sense (see the next section). This fact,
however, is somewhat hidden in the present proof, because a ``shortcut''
is used to obtain more rapidly this convexity result. The reader is
refered to \cite{efj}, where a special case (when the ambiant space has
constant curvature) is proved using an analog of theorem \ref{thpr}. The
resulting arguments are longer and more technical, but perhaps more
illuminating than those given here.

\bigskip

It is not clear whether the hypothesis concerning the gradients of the
curvature are really necessary here. On the other hand, the inequlities
on $K_1, K_2$ and $K_3$ are more or less optimal, as is pointed out in
section 8 using some examples. 

\bigskip

Note that a nice analog of theorem \ref{efimov} has been given by Smyth
and Xavier \cite{SX} in higher dimension, for hypersurfaces with Ricci
curvature conditions in $\R^{n+1}$ ($n\geq 3$). Some related results
have also been 
given by Smyth \cite{Smyth} in $S^{n+1}$. The approach they use,
however, is very different from the path followed here -- and it does
not seem to work at all for surfaces. It would be most interesting to
know whether something like the results of \cite{SX} applies to
hypersurfaces in Riemannian manifolds.

\section{How the proof works}

The proof of theorem \ref{thm-a} happens almost entirely on $\Si$ with
its third fundamental form, along with a compatible connection $\Nat$
which is defined in section 2. $\Nat$ is the Levi-Civita connection of
$\III$ when $M$ has constant curvature, but in general it has non-zero
torsion. Its torsion, however, is bounded. Section 2 contains the proof
of the following lemma, describing the basic geometric properties of
$\Nat$. 

\blm \label{Sit}
Under the hypothesis of theorem \ref{thm-a}, $\Nat$ is compatible with
$\III$, and has torsion $\tau$ bounded above by a constant
$\tau_0$. Its curvature $\Kt$ is bounded between two positive constant:
$$ K_5\geq \Kt\geq K_4>0~. $$
Moreover:
$$ 4K_4>\tau_0^2 ~. $$
\elm

We will also use the asymptotic directions of the immersion. More
precisely, we can suppose that $\Si$ is simply connected (otherwise
consider its universal cover, which again has an isometric immersion
into $M$). Therefore, we can choose two vector fields $U$ and $V$,
parallel to the asymptotic directions of the immersion, with unit norm
for $\III$. Since $U$ and $V$ are never parallel, we also demand that
$\angle(U,V)\in (0, \pi)$. Section 2 repeats this definition, and contains
the proof of the next lemma, about some key properties of $U$ and $V$.

\blm \label{UV}
There exists a constant $\tau_1>0$ such that the asymptotic vectors $U$
and $V$ satisfy:
$$ \|\Nat_U V\| \leq \tau_1 \sin(\angle(U,V)) \; , \; \; \|\Nat_V U\| \leq \tau_1\sin(\angle(U,V))~. $$
\elm

Section 3 contains some technical propositions concerning surfaces with
connections having bounded torsion. Section 4 is about an
amusing technical lemma which states that, if an asymptotic curve is
``almost closed'', then a propagation phenomenon happens. 
This is used in sections 5 and 6, which contain what is maybe
the central point of this paper. One must first define the convexity of
a (non-complete) surface in the following fairly natural way, basically
stating that a geodesic segment can not touch the boundary at an
interior point:

\begin{df} \label{df-cvx}
Let $(S, \dr S)$ be a surface, with a metric $g$ and a compatible
connection $D$. We say that $S$ is {\bf convex} if, when $(\ga_n)_{n\in \N}$
is a sequence 
of geodesic segments, $\ga_n:[0, L]\rightarrow \Si$, such that
$(\ga_n(t))$ converges in $\Sib$ for each $t\in [0,L]$, and when there
exists $t_0\in ]0, L[$ such that $\lim_{n\rightarrow \infty}
\ga_n(t_0)\in \dr \Si$, then 
$\lim_{n\rightarrow \infty} \ga_n(t)\in \dr \Sigma$ for all $t\in [0,
L]$.
\end{df}

Then:

\begin{lemma} \label{lem-e}
$\Sigma$, with $\III$ and $\Nat$, is convex.
\end{lemma}

To reach this goal, we define a specific notion of ``concavity'' of
$\dr_{\III}\Si$, and then prove that ``concave'' points are not
possible. The convexity of $\Si$ will then follow. First we choose
positive real numbers $k$ and $C$ and a point $x\in \dr_{\III}\Si$.

\begin{df} \label{df-nd}
A {\bf $(k,C)$-concave map} at $x$ is a map $\phi:[-d,d]\times
[0,d]\rightarrow \Sib$, with $d>0$, such that: 
\begin{itemize}
\item $\phi(0,0)=x$, and $\phi([-d,d]\times [0,d]\setminus (0,0))\subset
\Si$, $\phi$ being a smooth diffeomorphism on its image outside $(0,0)$; 
\item for each $y\in (0,d]$, the curve $\phi([-d,d]\times \{ y\})$ has
geodesic curvature $\kappa$ between $k$ and $Ck$, with its convex side
towards $x$, and $|\dr_2\kappa|\leq C$;
\item at each point of $[-d,d]\times [0,d]\setminus (0,0)$, $\dr_1\phi$
is orthogonal to $\dr_2\phi$, and $1\leq\|\dr_1\phi\|, \|\dr_2\phi\|\leq
C$. 
\end{itemize}
$d$ is called the {\bf diameter} of $\phi$ and is written as
$\diam(\phi)$. 
\end{df}

\begin{df} \label{df-conc}
Let $x_0\in \dr_{\III}\Si$. $\Si$ is {\bf $(k,C)$-concave} at $x_0$ if there
exists a $(k,C)$-concave map $\phi$ at $x_0$. $\Si$ is
{\bf concave} at $x_0$ if it is $(k,C)$-concave for some $k>0$ and $C>0$.
\end{df}

The point of this definition is the following result, which is proved in
section 5 using a technical lemma from section 4:

\blm \label{dom}
Under the hypothesis of theorem \ref{thm-a}, $\dr_{\III}\Si$ has no
concave point. 
\elm

On the other hand, it is proved in section 6 that:

\blm \label{ccv-cvx}
If $\Si$ has no concave point, then it is convex.
\elm

The proof of lemma \ref{lem-e} clearly follows from those two lemmas.
It is then proved in section 7 that:

\blm \label{borne-aire}
Under the hypothesis of theorem \ref{thm-a}, if $(\Si, \III, \Nat)$ is
convex, then it has bounded area. 
\elm

A contradiction will follow, because, by the Gauss formula, the ratio of
the area elements on $\Si$ for $I$ and for $\III$ is equal to the
absolute value of the
extrinsic curvature of the immersion, which is supposed to be bounded
away from $0$ in theorem \ref{thm-a}; and the area of $(\Si, I)$ is
infinite because $(\Si, I)$ is complete, simply connected, and with
negative curvature. 

\bigskip

{\bf Conventions:} in the whole paper, if $c:[a,b]\rightarrow \Si$ is a
piecewise smooth curve, and if $W\in T_{c(a)}\Si$, we let $\Pi(c;
W)$ be the parallel transport of $W$ at $c(b)$ along $c$.
Unless otherwise stated, all curves are parametrized at unit speed.

\section{Isometric immersions of surfaces}

This section contains some elementary results concerning the objects
induced on a $\Si$ by an immersion in a Riemannian
3-space $M$. We call $I$ the induced metric, $\Na$ its Levi-Civita
connection, and $\Na^M$ that of $M$.  

We suppose that $\Si$ is
contractible and oriented -- otherwise, consider its universal cover.
We can therefore choose a unit normal vector field $N$ to $\Si$, and
define a bundle morphism (the ``shape operator''):
\begin{eqnarray}
B: & T\Si \rightarrow & T\Si \nonumber \\
 & x \mapsto & \Na_x^M N~. \nonumber
\end{eqnarray}
It easy to check that $B$ is symmetric. From there follows the
definition of the third fundamental form of the immersion:
$$
\forall s\in \Si, \; \forall x,y \in T_s \Si, \; \III(x,y)=I(Bx,By)~. 
$$
If $M=\R^3$, then $\III$ is the pull-back of the canonical metric on
$S^2$ by the Gauss map. 

Let $R$ be the Riemann curvature tensor of $M$. Then $B$ satisfies the
following classical equations (see \cite{GHL} or 
\cite{S}, vol. III):
\beq
\forall s\in \Si, \; \forall x,y \in T_s\Si, \; (d^{\Na} B)(x,y) = 
-R_{x,y}n \label{CM}~,
\eeq
which is known as the Codazzi-Mainardi equation, and the Gauss equation:
\beq 
\forall s\in \Si, \;
\det(B_{s})= K_{e}:=K(s) - K_{M}(T_s\Si)~,
\label{Gauss}  
\eeq
where $K(s)$ is the curvature of $\Na$ at $s$.

The main point of this section is that the immersion also defines
on $\Si$ a connection which is compatible with $\III$, but in general
has torsion. 

\bdf
Let $\Nat$ be the connection defined on $\Si$ by:
$$
\Nat_{x}y = B^{-1} \Na_{x}(By)~.
$$
\edf

Remember that the torsion of a connection is a 2-form with value in the
tangent space, which is defined as:
$$ \tau(x,y):=\Nat_xy-\Nat_yx-[x,y] $$
The Levi-Civita connection of a metric is defined as the only compatible
connection with zero torsion. Note that, for Riemannian surfaces,
2-forms can be identified with functions, so we will often here consider
the torsion $\tau$ as a vector field on $\Si$. This identification will
always be made using, as a Riemannian metric, the third fundamental form
$\III$. 

The main property of $\Nat$ is given in the following proposition.

\bprop
$\Nat$ is compatible with $\III$. Its torsion
is bounded by:
$$
\forall s\in \Si, \; \forall x,y\in T_{s}\Si, \; \| \tau(x,y) \|_{\III}
\leq \| (d^{\Na}B) (x,y)\|_{I}~.
$$
\eprop

\bpv
Let $x,y$ and $z$ be three vector fields on $\Si$. Then:
\begin{eqnarray}
x.\III(y,z) & = & x.I(By, Bz) \nonumber \\
& = & I(\Na_x(By), Bz) + I(By, \Na_z(Bz)) \nonumber \\
& = & \III(B^{-1}\Na_x(By), z) + \III(y, B^{-1}\Na_x(Bz)) \nonumber \\ 
& = & \III(\Nat_x y, z) + \III(y, \Nat_x z)~, \nonumber
\end{eqnarray}
so that $\Nat$ is compatible with $\III$.

From the definition of $\Nat$:
\begin{eqnarray}
\tau(X,Y) & = & \Nat_{X}Y -\Nat_{Y}X -[X,Y] \nonumber 
\\
& = & B^{-1}\Na_{X}(BY) - B^{-1}\Na_{Y}(BX) - [X,Y] \nonumber \\
& = & B^{-1} (d^{\Na}B)(X,Y)~, \nonumber 
\end{eqnarray}
and this shows that:
\begin{eqnarray}
\|\tau\|_{\III}^{2} & = & \III(B^{-1} (d^{\Na}B)(X,Y),B^{-1} \nonumber 
(d^{\Na}B)(X,Y))  \nonumber \\
& = & I((d^{\Na}B)(X,Y),(d^{\Na}B)(X,Y))  \nonumber \\
& = & \|(d^{\Na}B)(X,Y)\|_{I}^{2}~, \nonumber 
\end{eqnarray}
and the result follows.
\epv

As a consequence, $\Nat$ is the Levi-Civita connection of $\III$ when
$M$ has constant curvature. When $M$ does not have constant curvature,
the previous proposition leads to the following control on the torsion
of $\Nat$:

\bprop \label{pr-tau}
If $(\Si, I)$ has curvature $K\leq K_1$, then, at any point $s\in \Si$,
the   torsion of $\Nat$ is bounded by:
\beq \label{etau0}
\|\tau\|_{\III} \leq \tau_0(\phi(s)) := \frac{K_M-K_m}
{2\sqrt{(K_m-K_1)(K_M-K_1)}}~, 
\eeq
where $K_m$ and $K_M$ are the minimum and the
maximum of the sectional curvatures of $M$ on tangent 2-planes at
$\phi(s)$. 
\eprop

\bpv
Let $(e_{1},e_{2})$ be the orthonormal basis of $T_{s}\Si$ for $I$ which
diagonalizes $B$, and let $k_{1},k_{2}$ be the associated
eigenvalues. We need to prove the upper bound above with $\tau$
replaced by $\tau(e_{1},e_{2})/(k_{1}k_{2})$, because $\tau$ is
skew-adjoint, and $((1/k_{1})e_{1},
(1/k_{2})e_{2})$ is an orthonormal basis of $T_{s}\Si$ for $\III$.

According to the previous proposition and to (\ref{CM}),
it is enough to prove that, under our curvature assumptions, for any
$m\in M$ and for any orthonormal basis $(x,y,n)$ of
$T_mM$:
$$ \frac{\| R(x,y)n \|_{I}}{K(x,y)-K_1} \leq \tau_0~, $$
where $K(x,y)$ is the sectional curvature of $M$ on the 2-plane
generated by $x$ and $y$.
Let $\Rc: \La^2M \rightarrow \La^2M$ the curvature operator, and $\mu$
the metric on $\La^2M$ coming from the metric on $M$. 
We need to prove that, for any $m\in M$, if, when $v,w \in \La^2_{m}M$
are orthogonal and have unit norm, $K_m \leq \mu(\Rc v,v) \leq K_M$,
then, with the same hypothesis on $v$ and $w$, we have:
$$ \left|\frac{\mu(\Rc v,w)}{\mu(\Rc v,v)-K_1} \right| \leq
\tau_0(m)~. $$ 

Let $m\in M$, and let $P\subset \La^2_mM$ be a 2-plane. Denote by
$Q$ the restriction of  $\Rc$ to $P$ followed by the orthonormal
projection on $P$, $p_1,p_2$ its eigenvectors, and  
$q_1,q_2$ its eigenvalues. If $v,w\in P$
are orthogonal with unit norm, they can be written as $v=\cos(\theta)p_1+
\sin(\theta)p_2$ and $w=\sin(\theta)p_1-\cos(\theta)p_2$, so that:
$$ \frac{\mu(Qv,w)}{\mu(Qv,v)-K_1} = \frac{(q_1-q_2)\cos(\theta)
\sin(\theta)}{q_1\cos^2(\theta)+q_2\sin^2(\theta)-K_1}~. $$
If now $\alpha:=\cos^2(\theta)$, we find that:
\beq \left|\frac{\mu(Qv,w)}{\mu(Qv,v)-K_1}\right|^2 = \frac{(q_1-q_2)^2 
\alpha(1-\alpha)}{
((q_1-q_2)\alpha+q_2-K_1)^2} \label{maxi}~. \eeq
This is maximal when:
$$ \alpha=\frac{K_1-q_2}{2K_1-q_1-q_2} $$
(which is in $[0,1]$ and corresponds to a possible value of
$\cos(\theta)$). Replacing $\alpha$ by this value in
(\ref{maxi}) shows that:
$$ \left|\frac{\mu(Qv,w)}{\mu(Qv,v)-K_1}\right| \leq
\frac{|q_1-q_2|}{2\sqrt{(q_1-K_1) 
(q_2-K_1)}}~. $$
Since the right side is maximal for $\{q_{1},q_{2}\}=\{K_{m},
K_{M}\}$, we find the upper bound we need for
$$ \frac{\mu(Qv,w)}{\mu(Qv,v)-K_1}~, $$
and the result for $\III(\tau, \tau)$ follows.
\epv

The previous proposition gives us informations about $\Nat$. Call $K_I$
the curvature of $I$ and $K_e$ the extrinsic curvature of the
immersion. Then:

\bcr \label{cor-Kt}
$\Nat$ is a connection compatible with $\III$, its torsion $\tau$
is bounded (for $\III$), at $s\in S$, by $\tau_0(\phi(s))$ (where
$\tau_0$ comes from (\ref{etau0})), and its curvature is:
\beq
\Kt = \frac{K_I}{K_e}~, \nonumber 
\eeq
with:
$$ 0 < K_4 \leq \Kt \leq K_5~, $$
where:
\begin{eqnarray}
K_5 = 1 \; \; \mbox{if}\;  \; K_2 \geq 0 & , & K_5 = \frac{K_1}{K_1-K_2}
\; \; 
\mbox{if} \; \; K_2 \leq 0  \nonumber \\
K_4 = 1 \; \; \mbox{if} \; \; K_3 \leq 0 & , & K_4 = \frac{K_1}{K_1-K_3}
\; \; 
\mbox{if} \; \; K_3 \geq 0~.  \nonumber 
\end{eqnarray}
\ecr

\bpv
We only have to prove the second assertion, concerning the curvature. 
Let $dv_I$ and $dv_{\III}$ be the area elements associated to the
metrics $I$ and $\III$ on $\Si$. By the Gauss formula
(\ref{Gauss}):
$$ dv_{\III} = K_e dv_I~. $$

Let $(e_1, e_2)$ be an orthonormal moving frame on $(\Si, I)$, and let
$\om$ be its connection 1-form, that is:
$$ \om(x):=I(\Na_x e_1, e_2)=-I(\Na_x e_2, e_1) $$
Then:
$$ K_I dv_I = \Om_I = -d\om~. $$
But $(B^{-1}e_1, B^{-1}e_2)$ is an orthonormal moving frame on $(\Si,
\III)$, and its connection 1-form $\om_{\III}$ is:
$$ \om_{\III}(x) = \III(\Nat_x (B^{-1}e_1, B^{-1}e_2)=\om(x)~. $$
Therefore:
$$ \Kt dv_{\III} = \Om_{\III} = -d\om_{\III} = -d\om = K dv_I~. $$
Those equations give the relation we need between
$\Kt$, $K_{e}$ and $K_I$.

The inequalities on $\Kt$ are direct consequences of this formula,
because: 
$$ \Kt = \frac{K_{I}}{K_e} \leq \frac{K_{I}}{K_{I}-K_2} $$
Now the function: $x\mapsto x/(x-\alpha)$ has as derivative: $x\mapsto
-\alpha/(x-\alpha)^2$, so its increasing for $\alpha\leq 0$ and
decreasing for $\alpha\geq 0$; for $\alpha=K_2$
we find the upper bound on $\Kt$ is obtained:
\bit
\item if $K_{2}\leq 0$, when $K_{I}\rightarrow K_{1}$, and it is
$K_{1}/(K_{1}-K_{2})$~; 
\item if $K_{2}\geq 0$, when $K_{I}\rightarrow \infty$, and it is $1$.
\eit
The same argument gives the lower bound for $\Kt$, with
$K_{2}$ remplaced by $K_{3}$.
\epv

Lemma \ref{Sit} is a direct consequence of proposition \ref{pr-tau} and
corollary \ref{cor-Kt}. 

\bigskip

We will now give two simple results which will be useful in the
sequel. First, $\Bt:=B^{-1}$ satisfies on $(\Si,\III,
\Nat)$ an equation similar to that
satisfied by $B$ on $(\Si, I)$ but even simpler: 

\bprop \label{pr-CM*}
On $(\Si,\III)$:
\beq
d^{\Nat}\Bt = 0 \label{CM*}
\eeq
\eprop

\bpv
A direct computation shows that, for $s\in \Si$ and $X,Y\in T_{s}\Si$~:
\begin{eqnarray}
(d^{\Nat}\Bt)(X,Y) & = & \Nat_{X}(\Bt Y) - \Nat_{Y}(\Bt X)
-B^{-1}([X,Y])  \nonumber \\
& = & B^{-1} \Na_{X}(B B^{-1} Y) - B^{-}\Na_{Y}(B B^{-1} X) - B^{-1} [X,Y]
\nonumber \\ 
& = & B^{-1} 0~, \nonumber 
\end{eqnarray}
because $\Na$ is torsion-free.
\epv

\bigskip

We will now describe some properties of $\Bt$ which will be useful later
on. Remember that, since $\det(\Bt)<0$, there exist at each point of
$\Si$ two vectors $U,V$ which have unit norm for $\III$, and such that:
\begin{eqnarray}
\Bt U & = & k J_{\III}U \nonumber 
\\
\Bt V & = & -k J_{\III}V~, \nonumber
\end{eqnarray}
where $J_{\III}$ is the complex structure defined by $\III$, and:
$$ 
k=|\det(\Bt)|^{1/2}=(K_{M}-K_{\Si})^{-1/2} ~.
$$ 
$U$ and $V$ are a priori defined only up to their orientation, but,
since we have supposed that $\Si$ is contractible, we can decide that,
in the remaining of this paper, $U$ et $V$ will be two globally defined
vector fields, oriented so that $\angle(U,V)\in ]0,\pi[$. 

\begin{remark} \label{UV-I}
The norms of $U$ and $V$ for $I$ are at most $1/\sqrt{K_2-K_1}$.
\end{remark}

\bpv
By definition:
$$
I(U,U)=\III(\Bt U, \Bt U)=\III(kJ_{\III}U,
kJ_{\III}U)=k^2=(K_M-K_{\Si})^{-1} \leq (K_2-K_1)^{-1}.
$$
\epv

\bigskip

The remainder of this section is dedicated to some elementary facts
about the asymptotic curves of the immersion, as seen on $(\Si, \III,
\Nat)$. Those curves have been well studied on $(\Si, I)$; for instance,
they have been used before \cite{Ef1} in \cite{Ef3} to prove that there
exists a constant $k$ such that, if a smooth, complete Riemannian
surface $S$ has 
uniformly negative curvature, and if the norm of the gradient of this
curvature is bounded by $k$, then $S$ has no isometric immersion into
$\R^3$. But we only give here some details on the local behavior of
asymptotic curves on $(\Si, \III, \Nat)$.

In all this paper, $\theta$ denotes the angle between $U$ and $V$ for
$\III$. As above, we suppose that 
$\theta\in ]0,\pi[$. Note that $\theta$ is close to $0$ (or to $\pi$)
when the immersion $\phi$ is ``degenerate'': the mean curvature of
$\phi$ is $\cot(\theta)(|\det(\Bt)|)^{-1/2}$.

\bprop
At each point of $\Si$:
\begin{eqnarray}
\Nat_{V} U & = & - \frac{\sin(\theta)}{2} (U.\kappa +
\III(\tau,J_{\III}U)) J_{\III}U \label{UV*} \\
\Nat_{U} V & = & \frac{\sin(\theta)}{2} (V.\kappa + 
\III(\tau,J_{\III}V)) J_{\III}V~, \label{VU*} 
\end{eqnarray}
with $\kappa=\ln(k^{-1})=-\ln(k)$
\eprop

\bpv
From (\ref{CM*}):
$$ (d^{\Nat} \Bt)(U,V) = 0~, $$
so, if 
$\om_U:=\III(\Nat_{U}V, J_{\III}V)$ 
and $\om_V:=\III(\Nat_{V}U, J_{\III}U)$:
$$ \Nat_{U}(\Bt V)-\Nat_{V}(\Bt U)-\Bt(\Nat_{U}V-\Nat_{V}U
- \sin(\theta)\tau) = 0~, $$ 
so that:
$$ - \Nat_U (k J_{\III} V) -
\Nat_{V}(k J_{\III} U) -\Bt(\om_{U} J_{\III} V) +
\Bt(\om_{V} J_{\III} U) + \sin(\theta)\Bt\tau = 0~. $$
But $\sin (\theta)J_{\III} U=V - \cos (\theta)U $
and $\sin(\theta)J_{\III} V=\cos (\theta)V - U$, and it follows that:
$$ \om_V U + \om_U V + \left(- V.\kappa +\frac{\om_U}{\sin(\theta)}
-\frac{\om_V \cos (\theta)}{\sin (\theta)}\right)J_{\III} U+ \coupeq
+ \left(- U.\kappa +\frac{\om_U \cos (\theta)}{\sin
(\theta)}-\frac{\om_V}{\sin (\theta)}\right) J_{\III} V +
\frac{\sin(\theta)}{k} \Bt\tau = 0~. $$ 
Take the scalar product (for $\III$)
with $U$ and then with $V$, and use the symmetry of $\Bt$ with respect
to $\III$ to obtain the result.
\epv

We will use this proposition to show that $U$ and $V$ each behave well
along the integral curves of the other. This will be used in section 4
to obtain a key technical lemma on asymptotic curves. Note that the
hypothesis of theorem \ref{thm-a} on the gradient of the curvature
appears only here. 

Remember that, according to the hypothesis of theorem \ref{thm-a}:
\bit
\item There exists $c_{\si}>0$ such that, for all $s\in \Si$ and all
$x\in T_{s}\Si$:
$$
\|x.K_{\si}\| \leq c_{\si} \|x\|_{\si} |K_{\si}|^{3/2}~.
$$ 
\item There exists $c_{\mu}>0$ such that, for all $m\in M$ and all $x\in
T_{m}M$, for each 2-plane $P\in G^{2}_{m}M$:
$$
|(\Na^{M}_{x}K_{\mu})(P)|\leq c_{\mu}\|x\|~.
$$ 
\eit

Then:

\bcr \label{cr-UV}
There exists $\tau_{1}>0$ (depending on $K_{1},K_{2},K_{3}, c_{\si},
c_{\mu}$ only)  such that:
\begin{eqnarray}
\|\Nat_{U}V\|_{\III} & \leq & \tau_{1} |\sin(\theta)| \label{UV2*}
\\
\|\Nat_{V}U\|_{\III} & \leq & \tau_{1} |\sin(\theta)|~. \label{VU2*}
\end{eqnarray}
\ecr

\bpv
According to the previous proposition:
\begin{eqnarray}
\| \Nat_{U}V \|_{\III} & \leq & \left|\frac{\sin(\theta)}{2}\right|
\left( \left| \frac{V.K_{e}}{2K_{e}} \right| + \|\tau \|_{\III}
\right) \nonumber  
\\
& \leq & \left|\frac{\sin(\theta)}{4}\right| \left( \left|
\frac{V.K_{\si}}{K_{e}}\right| + \left| \frac{V.K_{\mu}}{K_{e}}
\right| + 2 \tau_{0} \right)~.  \nonumber 
\end{eqnarray}
Let $x\in M$, call $K_\mu^x$ the restriction of $K_\mu$ to the Grassmannian
of 2-planes in $T_xM$. A simple compactness argument shows that there
exists a constant $C_M$ (which does not depend on $M$) such that:
$$ dK_\mu^x \leq C_M K_M^x~, $$
where $K_M^x$ is the maximum of the sectional curvatures of $M$ at $x$.
Therefore, isolating in $V.K_{\mu}$ a part coming from the derivative of
$K_{\mu}$ 
from another coming from the rotation of the tangent plane during a
displacement in the direction of $V$ shows that:
$$
\| \Nat_{U}V \|_{\III} \leq \left|\frac{\sin(\theta)}{4}\right| \left(
\left| \frac{V.K_{\si}}{K_{e}} \right|+ \left|
\frac{(\Nat_{\phi_{*}V}K_{\mu})(\phi_{*}(T_{s}\Si))}{K_{e}}\right| +
\|V\|_{\III} C_M \left|\frac{K_{M}}{K_{e}}\right| + 2 \tau_{0}
\right)~,
$$
because the norm of the rotation of $\phi_{*}T_{s}\Si$ during
displacements along $\Si$ is measured by $\III$. But
$\|V\|_{\III}=1$ and $\|V\|_{I}=k=K_e^{-1/2}$, so:
$$
\| \Nat_{U} V\|_{\III}\leq \left|\frac{\sin(\theta)}{4}\right| \left(
k^2 |V.K_{\si}| + 
k^2 |(\Nat_{\phi_{*}V} K_{\mu})(\phi_{*}(T_{s}\Si)) + \left|
\frac{C_M K_{M}}{K_{e}} \right| + 2 \tau_{0} \right)~,
$$
and, if $k_{M}$ is the maximal possible value of $k$, i.e.
$k_{M}=(K_{2}-K_{1}^{-1/2}$:
$$
\| \Nat_{U} V\|_{\III}\leq \left|\frac{\sin(\theta)}{4}\right| \left(
k_M^3 c_{\si}  + k_M^3 c_{\mu} +
k_M^2 C_M |K_{3}| + 2 \tau_{0} \right)~,
$$
whence the first result. The same computation with $U$ and $V$
interchanged gives the same bound for $\| \Nat_{V} U \|_{\III}$.
\epv

Lemma \ref{UV} is no more than a restatement of corollary \ref{cr-UV}.

\section{Connections with bounded torsion}

This section contains some simple technical propositions describing some
properties of surfaces with metrics and compatible connections with
bounded torsion.

First note that the Gauss-Bonnet theorem remains valid in this setting:
if $D$ is a compact, simply connected domain in $\Si$ with 
smooth boundary, then 
the integral of the geodesic curvature (for $\Nat$) of $\dr D$ is equal
to $2\pi$ 
minus the integral of the curvature $\Kt$ of $\Nat$ over $D$.

This is proved as follows. Let $(X,Y)$ be an orthogonal moving frame on $D\setminus
\{p\}$, where $p$ is a point in $D$, with $X$ tangent to $\dr D$ and to
the ``circles'' $\dr B(p, \eps)$ for $\eps$ small enough. Let $\om$
the connection 1-form of $(X,Y)$, and $\Om$ its curvature 2-form. By
definition of $\Kt$:
$$
\Om = \Kt dv~,
$$
where $dv$ is the area form of $\III$; moreover:
$$
\Om = -d \om~,
$$
so
$$
\int_{D} \Om = - \int_{\dr M}\om -\lim_{\eps\rightarrow 0} \int_{\dr
B(p,\eps)} \om~.
$$
Therefore, if $\kappa$ is the geodesic curvature of $\dr D$:
$$
\int_{D} \Kt dv = - \int_{\dr D} \kappa ds + 2\pi~.
$$

This theorem of course remains true if $\dr D$ is only piecewise
smooth, with the adequate contributions from the singular points. 

\bigskip

We now describe some properties of geodesics which ressemble
those for Jacobi fields along geodesics when the connection has no
torsion. But the torsion comes into the equations so that the usual
equalities are replaced by inequalities.

Let $(g_{s})_{s\in [0,1]}$ be a family of $\Nat$-geodesic,
$g_{s}:[0,L]\rightarrow \Si$, parametrized at unit speed. For each $s\in
[0,1]$ and each $t\in [0, 
L]$, we let $g':=\dr g_s(t)/\dr t$ and $\gd:=\dr g_s(t)/\dr s$. For
$s=0$, $\gd$ is a kind of Jacobi field along $g_0$, and we can call $x$
and $y$ the functions from $[0,L]$ to $\R$ such that, for $s=0$:
$$ \gd = x g' + y J_{\III} g'~. $$
We also call $\tau_x(t):=\III(\tau, g'_s(t))$ and $\tau_y:=\III(\tau,
J_{\III}g'_s(t))$. 

\bprop \label{edo-jac}
$x$ and $y$ are solutions of:
\beq \label{edo-x}
x' = y\tau_x 
\eeq
\beq \label{edo-y}
y'' = -\Kt y + (y\tau_y)'~. 
\eeq
\eprop

\bpv
By definition of $g'$ and $\gd$, $[g', \gd]=0$, so that, by definition
of the torsion:
\beq \label{int-1}
\Nat_{\gd} g' = \Nat_{g'}\gd - \tau(g', \gd)~. 
\eeq
Taking the scalar product with $g'$ and using the fact that the $(g_s)$
are parametrized at unit speed shows that:
$$ 0 = \gd.\III(g',g') = 2 \III(\Nat_{g'}\gd - \tau(g', \gd), g')~. $$
Therefore:
$$ g'.\III(\gd, g') - \III(\tau(g', \gd), g') = 0 $$
and we obtain the first equation.

Coming back to equation (\ref{int-1}), we see that:
\begin{eqnarray}
\Nat_{g'}\Nat_{g'}\gd & = & \Nat_{g'}\Nat_{\gd} g' + \Nat_{g'}(\tau(g',
\gd)) \nonumber \\
& = & R_{g',\gd} g' + \Nat_{\gd}\Nat_{g'} g' + \Nat_{g'}(\tau(g',
\gd)) \nonumber \\
& = & - \Kt y J_{\III} g' + \Nat_{g'}(y \tau_x g' + y \tau_y J_{\III}g')
\nonumber \\
& = & (y \tau_x)' g' + (-\Kt y + (y \tau_y)') J_{\III}g'~,
\nonumber 
\end{eqnarray}
and the second equations follows (as well as the derivative of the
first). 
\epv

\bcr \label{cr-ido}
There exists $t_g>0$, depending on $K_4, K_5$ and $\tau_0$, such that,
if $x(0)=y(0)=0$, then, for all $t\in [0, t_g]$:
\beq \label{ido-y}
\frac{y'(0) t}{2} \leq y(t) \leq 2y'(0) t
\eeq
\beq \label{ido-x}
|x(t)|\leq \tau_0 y'(0) t^2~.
\eeq
\ecr

\bpv
Integrating (\ref{edo-y}) shows that, for $t\in [0, L]$:
$$ y'(t)-y'(0) = \int_0^t -K(s) y(s) ds + y(t)\tau_y(t)~, $$
so that:
$$ -\tau_0 y(t) - K_5\int_0^t y(s)ds \leq y'(t) - y'(0)\leq \tau_0 y(t) -
K_4\int_0^t y(s)ds~. $$
Let: 
$$ t_1:=\inf\{ t\geq 0 \; | \; y(t)\not\in [y'(0)t/2, 2y'(0)t] \}~. $$
For $t\leq t_1$:
$$ y'(0)(1-2\tau_0 t-K_5t^2) \leq y'(t)\leq y'(0)(1+2\tau_0 t-K_4
t^2/4)~. $$ 
Thus there exists $t_g>0$ such that, if $t_1<t_g$, then:
$$ \frac{y'(0)t}{2} \leq y(t) \leq 2 y'(0) t~, $$
which contradicts the definition of $t_1$. So $t_1\geq t_g$, and
equation (\ref{ido-y}) follows. (\ref{ido-x}) is a direct
consequence using (\ref{edo-x}). 
\epv

\bcr \label{diff-loc}
If $x\in \Si$ and $v\in T_x\Si$ is a vector of norm at most $t_g$ at
which the exponential at $x$ for $\Nat$, $\exp_x^{\Nat}$, is defined,
then $\exp_x^{\Nat}$ is a local diffeomorphism at $v$.
\ecr

\bpv
Let:
$$ v':= \Pi(\exp_x^{\Nat}([0,1]v), v)~. $$
Equation (\ref{ido-y}) shows that:
$$ \III((d_v\exp_x^{\Nat})(J_{\III}v), J_{\III}v') \neq 0~, $$
while it is easy to check that:
$$ (d_v\exp_x^{\Nat})(v) = v'~, $$
because this corresponds to a change in the parametrization of the
geodesic starting at $x$ in the direction of $v$.
\epv

\bcr \label{ex-geod}
Let $\Om\subset \Si$ be an open subset with locally convex boundary,
$\Omb\subset \Si$. For any $x,y\in \Om$ with $d_{\III}(x,y)\leq t_g$,
there exists a unique $\Nat$-geodesic of length $d_{\III}(x,y)$ between
$x$ and $y$.
\ecr

\bpv
Let $\Om'$ be the inverse image of $\Om$ by the restriction of
$\exp_x^{\Nat}$ to the ball of radius $t_g$. By the previous corollary
and the local convexity of $\Om$, the restriction of $\exp_x^{\Nat}$ is
a diffeomorphism onto its image.
\epv

Here is another elementary corollary of proposition
\ref{edo-jac}.

\bcr \label{dfo-par}
For all $\eps>0$, there exists $\alpha>0$ such that , if $L\leq \alpha$
and $y'(0)=x(0)=0$, then:
\beq \label{par-y}
\forall t\in [0, L], |y(t)-y(0)|\leq \eps
\eeq
\beq \label{par-x}
\left| x(L) - y(0) \int_0^L \tau_y(s)ds \right| \leq \eps L
\eeq
\ecr

\bpv
(\ref{par-y}) is a simple consequence of (\ref{edo-y}), and (\ref{par-x})
then follows from (\ref{edo-x}).
\epv

We can now consider a family of geodesic rays starting from a
given point, and describe how they behave relative to one another. 
Let $(g_\theta)_{\theta\in [0,\theta_0]}$ be a family of maximal rays, with
$g_\theta:[0,L_\theta)\rightarrow \Si$, $L_\theta\in \R_+^*\cup
\{\infty\}$, and with $g_\theta(0)=g_0(0)$ and $\angle(g'_0(0),
g'_\theta(0))=\theta$ for each $\theta\in [0, \theta_0]$.

For $s\in [0, L_0)$, let $n_s$ be the maximal geodesic ray with
$n_s(0)=g_0(s)$ and $n'_s(0)=J_{\III}g'_0(s)$. Choose $s_1>0$ and
$\theta_1>0$, and suppose that there is no $\theta, t_0, s, u_0$ with
$s\leq s_1$ and $\theta\leq \theta_1$ such that:
$$ \lim_{t\rightarrow t_0}g_\theta(t) = \lim_{u\rightarrow u_0} n_s(u)
\in \dr_{\III}\Si~. $$
Then:

\bprop \label{fam-geod}
There exists a constant $S>0$ and, for each $\eps>0$ small enough and
each $s_1>0$, there exists $\Theta(\eps, s_1)>0$ (both also 
depending on $\tau_0, K_4, K_5$) such that, if $s\leq s_1$ and
$\theta\leq \Theta(\eps, s_1)$, then:
\begin{enumerate}
\item $g_\theta$ intersects $n_s$ at a point $n_s(u_\theta(s))$ (with
$g_\theta\cap n_s([0, u_\theta(s)))=\emptyset$);
\item the restriction of $|u_\theta|$ to
$[0, s]$ remains bounded by $\eps$;
\item if $s\geq S$, there exists $s'\in [0, s]$ such that
$u_\theta(s') = - \eps \theta$.
\end{enumerate}
\eprop

\bpv
Let $u_M$ be a small real number; we will see later how small $u_M$ has
to be. For $\theta\in [0, \theta_1]$, let:
$$ \alpha_\theta(s):=\angle(-J_{\III}n'_s(u_\theta(s)), g'_\theta)~. $$
Then define:
$$ s_\theta:=\sup\{s\in \R_+~ | ~ \forall s'\in [0, s],
|u_\theta(s')|\leq u_M ~ \mbox{and} ~ |\alpha_\theta(s')|\leq u_M\}~. $$
For $s\in [0,s_\theta]$, apply the Gauss-Bonnet theorem to an
infinitesimal strip bounded by $g_0([s, s+ds])$, $n_s([0,
u_\theta(s)])$, $n_{s+ds}([0, u_\theta(s+ds)])$ and $g_{\theta}$. This
shows that:
$$ \alpha'_\theta(s) = - \int_0^{u_\theta(s)} \Kt(n_s(t)) \left\|
n'_s(t)\wedge \frac{\dr}{\dr s} n_s(t) \right\| dt~, $$
so that:
\beq \label{eq-alp}
\alpha'_\theta(s) = -k(s) u_\theta(s)~, 
\eeq
where $k(s)\in [K_4-\eps, K_5+\eps]$ if $u_M$ is small enough (this last
step uses corollary \ref{dfo-par} applied to the family $(n_s)$).

Again by corollary \ref{dfo-par}, it is not hard to check that, again
for $u_M$ small enough:
$$ \left\| \frac{\dr}{\dr s}n_s(u_\theta(s)) +
J_{\III}n'_s(u_\theta(s))\right\| \leq \frac{\eps}{4}~. $$
Thus, with (\ref{par-x}):
$$
\left( 1-\frac{\eps}{4}\right) \sin\alpha_\theta(s) + (1-\eps)
\int_0^{u_\theta(s)} \tau(-J_{\III}n_s'(t)) dt \leq u'_\theta(s) \leq 
\coupeq \leq
\left( 1+\frac{\eps}{4}\right) \sin\alpha_\theta(s) + (1+\eps)
\int_0^{u_\theta(s)} \tau(-J_{\III}n_s'(t)) dt~.
$$
This can be written, for $u_M$ small enough, as:
\beq \label{eq-up}
u'_\theta(s) = \lambda(s)\alpha_\theta(s) + \tau(s) u_\theta(s)~, 
\eeq
with:
$$ |\lambda(s)-1|\leq \eps, ~~ |\tau(s)|\leq \tau_0 (1+\eps)~. $$
Let:
$$ X(s) := \left(
\begin{array}{c}
u_\theta(s) \\
\alpha_\theta(s)
\end{array}
\right)~.
$$
Then:
$$ X'(s) = m(s) X(s)~, $$
with:
$$ m(s) := \left(
\begin{array}{cc}
\tau(s) & \lambda(s) \\
-k(s) & 0
\end{array}
\right)~.
$$
Thus, by integration:
$$ X(s) = \exp(s M(s)) X(0)~, $$
where:
$$ M(s) := \left(
\begin{array}{cc}
T(s) & \Lambda(s) \\
-K(s) & 0
\end{array}
\right)~,
$$
with:
$$ |T(s)|\leq (1+\eps)\tau_0, ~ |\Lambda(s)-1|\leq \eps, ~ K_4-\eps
\leq K(s) \leq K_5+\eps~. $$

The eigenvalues of $M(s)$ are the roots of:
$$ X(X-T(s))+\Lambda(s) K(s) = 0~. $$
If $\eps$ is so small that $4(K_4-\eps)(1-\eps)>(1+\eps)^2\tau_0$, those
roots can be written as
$\alpha \pm i\beta$, where:
$$ |\alpha|=\frac{T(s)}{2} \leq \frac{(1+\eps)\tau_0}{2}, ~ |\beta|\geq
\frac{\sqrt{4(K_4-\eps)(1-\eps)-(1+\eps)^2\tau_0^{2}}}{2}~. $$

Therefore, in a well chosen frame, the orbits of $X(s)$ are ``spirals''
around $0$, with an angular speed which is bounded from below. This
already proves, with the upper bound on $\alpha$, that, if $\theta$ is
smaller than some $\Theta(\eps, s)$, 
then $s_\theta\geq s$, so that $u_M$ is not reached and the computations
above hold on all of $[0, s]$. This proves point (2).

Moreover, the trajectories $(X(s'))_{s'\in [0, s]}$ can not remain in a
half-plane, so that $u_\theta$ has to become negative after a time which
is bounded in term of $\beta$ (which itself is bounded from below). This
leads to point (3) of the proposition.
\epv

Finally, the same kind of argument will show the following similar
proposition, which deals with convex curves instead of geodesics. The
proof is similar to the one we have just finished, so it is described
somewhat faster.

\bprop \label{prop-gagat}
Let $S$ be a convex domain in $\Si$, with boundary $\dr S$ containing as
connected components two complete curves $\ga$ and $\gat$.
Suppose that $K_{4}>\tau_{0}^2/4$. Then $d_{\III}(\ga, \gat)>0$.
\eprop

\bpv
If $\ga$ or $\gat$ is compact, the result is
obvious, so we suppose here that neither $\ga$ nor $\gat$
is compact. The proof is by contradiction, so we suppose that
$d_{\III}(\ga, \gat) = 0$.

First note that a rather direct smoothing argument shows that, for any
$\eps_r>0$, there are smooth curves $\ga_r, \gat_r:\R\rightarrow \Si$
such that: 
\begin{itemize}
\item $(\dr S\setminus (\ga\cup \gat))\cup (\ga_r\cup \gat_r)$ bounds a
connected closed set $S_r$ which contains $S$;
\item for each $s\in \R$, the curvatures $\kappa(t)$ and $\kappat(t)$ of
$\ga_r$ at $\ga_r(t)$ and of $\gat_r$ at $\gat_r(t)$ respectively are
bounded by:
$$ \kappa(t)\geq -\eps d_{\III}(\ga_r(t), \gat_r)~, ~~
\kappat(t)\geq -\eps d_{\III}(\gat_r(t), \ga_r)~, $$
where both curvatures are with respect to the normal oriented towards
the interior of $S'$;
\item $\liminf_{t\rightarrow \infty} d_{\III}(\ga_r(t), \gat_r)=0$.
\end{itemize}

For $s\in \R$, let:
$$ d(s) = d_{\III}(\ga_r(s), \gat_r)~.  $$
Thus $d$ is not bounded away from $0$ near $+\infty$.

Choose $\eps>0$. There exists $s_0\in \R$ with
\beq \label{ddp} 
d(s_0) \leq \eps, \; d'(s_0) \leq \eps~.
\eeq
If $\eps$ is small enough, it is not difficult to show, using
\ref{ex-geod}, that there exists a $\Nat$-geodesic $n_{s_0}$ connecting 
$\ga_r(s_0)$ to $\gat_r$, of length at most $2\eps$, orthogonal to
$\gat_r$. For $s>s_0$, let $n_s$ be the maximal $\Nat$-geodesic
starting at $\ga_r(s)$ with speed equal to the parallel transport of
$n_s'(0)$ at $\ga_r(s)$ along $\ga_r$.

Let $r(s)$ be the distance along $n_s$ between
$\ga_r(s)$ and the first intersection of $n(s)$ with $\gat_r$, $\beta(s)$
the angle between $-J_{\III}n_s'(0)$ and $\ga_r'(s)$, $\alpha$ the angle
between $-J_{\III}n_s'(r(s))$ and $\gat_r'$. By construction,
$\alpha(s_0)=0$, while, by corollary \ref{dfo-par} and (\ref{ddp}),
$\beta(s_0)$ is small.
Let $u_M$ be again a small real number, for which precisions will come
later. Define:
$$ s_M:=\sup\{s\geq s_0 ~ | ~ \forall s'\in [s_0, s],
|u(s')|\leq u_M ~ \mbox{and} ~ |\alpha(s')|\leq u_M ~ \mbox{and} ~
|\beta(s')|\leq u_M\}~. $$ 
The definition of $\beta$ and the ``almost'' convexity of $\ga_r$ show
that $\beta'(s)\geq -2 \eps u(s)$, while the same
application of the Gauss-Bonnet theorem as the one leading to
(\ref{eq-alp}) shows again that $\alpha'(s)=-k(s) u(s)$, but with only
$k(s)\geq K_4-2\eps$, while the upper bound is lost because $\gat_r$ is
only ``almost convex'' instead of geodesic.

Moreover, the same argument as the one leading to (\ref{eq-up}) shows
that:
$$ u'(s) = \lambda(s)(\alpha(s)-\beta(s)) + \tau(s) u(s)~, $$
again with:
$$ |\lambda(s)-1|\leq \eps, ~~ |\tau(s)|\leq \tau_0 (1+\eps)~. $$

The rest of the proof can now be done just as in the proof of
proposition \ref{fam-geod}, with $\alpha_\theta$ replaced by
$\alpha-\beta$, to obtain that there exists $S>0$ (depending on $K_4$
and $\tau_0$) such that: 
\begin{itemize}
\item either there exists $s\in [s_0, s_0+S]$ such that $u(s)=0$, and
this proves 
the proposition;
\item or $s_M<s_0+S$, and in this case the upper bound on the norm of $X$
shows that, if $\eps$ has been chosen small enough, then either
$\alpha(s)=-u_M$ or $\beta(s)=u_M$.
\end{itemize}
But then, again for $\eps$ small enough, it is not difficult to show
that there exists $S'>0$ such that there exists $s\in (s_M, s_M+S')$
such that $u(s)=0$, so that the proposition holds also in that case.
\epv

\section{Asymptotic curves}

This section contains the proof of lemma \ref{propag}, a technical
statement which will have a central role later on. This lemma, along
with its proof, is 
similar to a lemma from \cite{efj}, but more detailed estimates are
necessary here. First, we introduce a
simple notation. It is written for an integral curve of $U$, but the
analog for an integral curve of $V$ should be obvious.

\bdf 
Let $\ga:[0, L]\rightarrow \Si$ be an integral curve of $U$ or
$V$. Then:
$$ \delta_\gamma:=\pi + \inf_{t\in [0,L]} \theta(\gamma(t)) - \sup_{t\in
[0,L]} \theta(\gamma(t))~, $$
and:
$$ \sigma_\gamma:=\int_0^L \sin(\theta(\gamma(s))) ds~. $$
\edf
Thus $\delta_\gamma\in (0, \pi)$; heuristically, because of (\ref{UV2*})
and (\ref{VU2*}), $\delta_\gamma$ is
small when $\ga$ has a segment which looks like a closed loop.
The following definition is very natural:

\bdf
Let $\eps>0$. A curve $c:[0,L]\rightarrow \Si$ is an
$\eps$-quasi-geodesic if, for each $s\in [0, L]$, the absolute value of
the angle between $c'(s)$ and 
$\Pi(c_{|[0,s]}; c'(0))$ is at most $\eps$.
\edf

\blm \label{propag}
There exists $T_0>0, C_0>0$ and $\eps_0>0$ as follows. 
Let $g$ be an integral curve of $U$ of length $L_{g}\leq T_0$,
with  $\eps:=\max(\delta_g, \sigma_g) \leq \eps_0$. Let 
$h_u:[-T_0,T_0]\rightarrow \Si$ be the integral curve of $V$ with
$h_u(0)=g(u)$. Then, for any $u\in [0, L_g]$, $h_u$ is a
$C_0\eps$-quasi-geodesic. 
\elm

\vspace{0.3cm}
\centerline{\psfig{figure=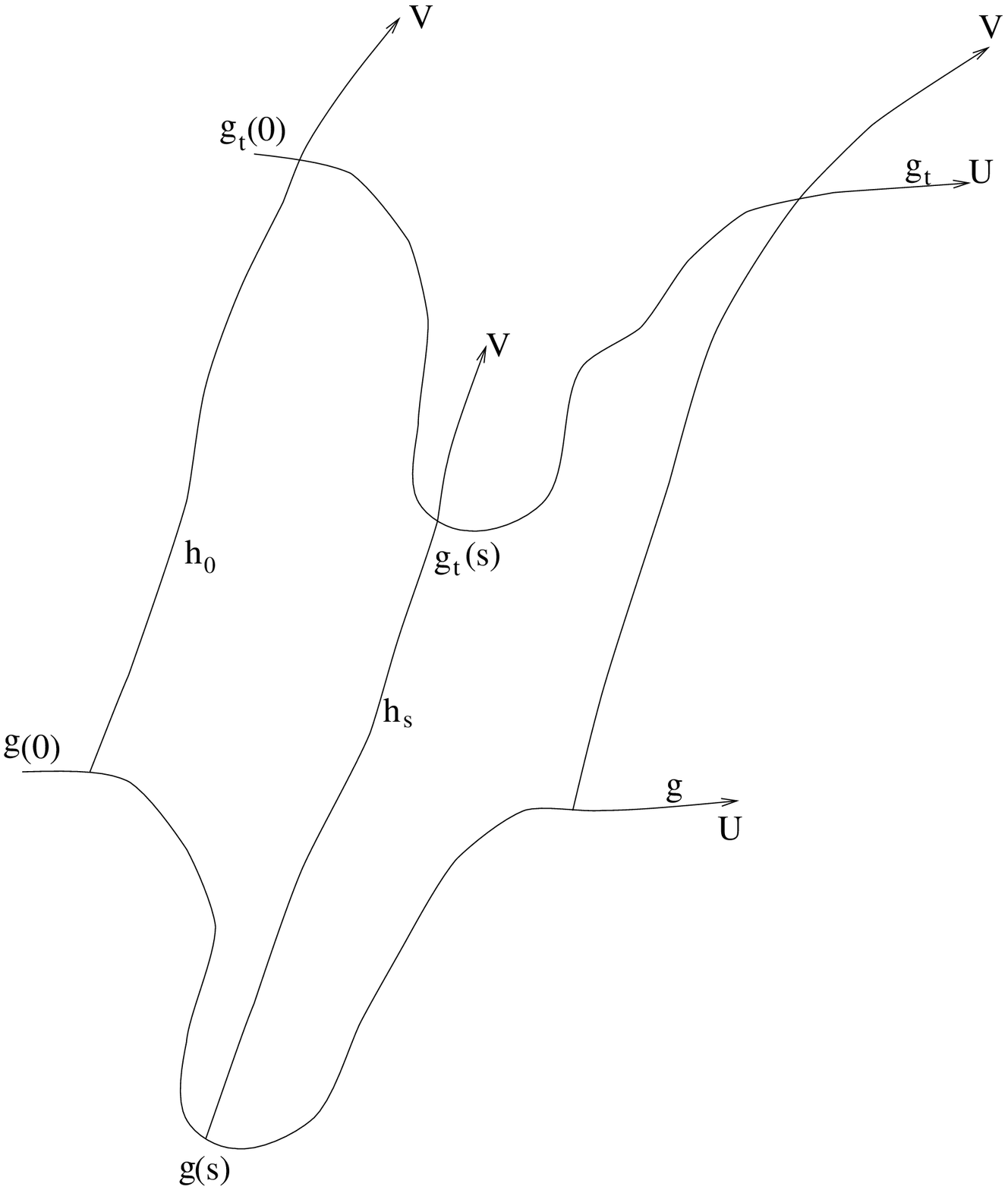,height=8cm}}
\centerline{\bf Figure 4.1} \vspace{0.3cm}

More could actually be said: the curve $g$ ``propagates'' along the flow
of $V$, 
that is, under this flow $V$, the integral curves of $U$ corresponding
(in some natural sense) to $g$ still have very small values of 
$\delta$ and of $\sigma$. This should be clear from the proof, although
we do not elaborate on it since it is not used later on.

The proof of this lemma rests on the following:

\bprop \label{controles}
There exists $L_M>0$ and continuous functions $\phi, \Phi: [0, L_M]\times [0,
L_M]\rightarrow \R_+$ such that, for all $L\in [0, L_M]$,
$\phi(0, L)=0$, and $\Phi(0, L)=\Phi(L,0)=0$, with the following
properties. 
If $g:
[0, L_0]\rightarrow \Si$ and $\gba: [0, L_1]\rightarrow \Si$ are integral
curves of $U$, if $h:
[0, L'_0]\rightarrow \Si$ and $\hb: [0, L'_1]\rightarrow \Si$ are
integral curves of $V$, with $g(0)=h(0)$, $g(L_0)=\hb(0)$,
$\gba(0)=h(L'_0)$ and $\gba(L_1)=\hb(L'_1)$, and if $L_0\leq L_M$ and
$L'_0\leq L_M$, then $L_1\leq \phi(L_0,
L'_0)$, $L'_1\leq \phi(L'_0, L_0)$, and the area of the domain bounded
by $g, \gba, h$ and $\hb$ is at most $\Phi(L_0, L'_0)$.
\eprop

\bpv
Let $u\in [0, L_0]$ and $v\in [0, L'_0]$. We can suppose that the integral
curve of $V$ starting at $g(u)$ meets the integral curve of $U$
starting from $h(v)$: otherwise, the proposition would fail slightly
before the first value of $u$ such that the intersection does not exist,
because then the length of both $\gba$ and $\hb$ would go to infinity.
We call $g_v(u)$ the intersection of the integral
curve of $V$ starting at $g(u)$ with the integral curve of $U$
starting from $h(v)$; this
intersection has to be unique because $\Si$ is simply connected and $U$
and $V$ are transverse. 

Let $\partial_u = \partial /\partial u, \partial_v = \partial
/\partial v$. Then:
$$
\dr_v g_v(u) = \alpha(u,v) V\; , \; \; \dr_u g_v(u) =
\beta(u,v) U~.
$$

By definition of $\tau$:
$$ \Nat_{\alpha V} (\beta U) - \Nat_{\beta U}
(\alpha V) - [\alpha V,\beta U] = \tau(\alpha V,\beta U)~, $$
so:
$$ (\dr_v \beta) U + \alpha \beta \Nat_VU - (\dr_u \alpha)V - \alpha
\beta \Nat_UV = - \alpha \beta \sin(\theta) \tau~. $$
Take the scalar product (for $\III$) of this equation with $J_{\III}U$
to obtain that:
$$
- \sin(\theta) \dr_u \alpha + \alpha \beta (\langle \Nat_{V} U, J_{\III}U
\rangle_{\III} - \langle \Nat_{U} V, J_{\III}U \rangle_{\III}) = - \alpha
\beta \sin(\theta) \langle \tau, J_{\III} U\rangle_{\III}~,
$$
which shows, along with lemmas \ref{Sit} and \ref{UV}, that:
$$
\left| \dr_u \alpha\right| \leq (\tau_0+2\tau_{1}) |\alpha \beta |~. 
$$ 
The same proof can be used to show also that:
$$
\left| \dr_v \beta \right| \leq  (\tau_0+2\tau_{1}) |\alpha \beta |~.
$$
In other terms:
\beq \label{UlVv}
|U.\alpha| \leq (\tau_0 + 2 \tau_1)\alpha ~ , ~ ~
|V.\beta| \leq (\tau_0 + 2 \tau_1)\beta~.
\eeq
Moreover, $\alpha(0, v)=1$ and $\beta(u, 0)=1$.

Integrate (\ref{UlVv}) over $g_v$ to obtain that:
$$ \alpha(u,v) \leq \exp ((\tau_0+2\tau_1) L(g_v))~. $$
Using (\ref{UlVv}) again leads to:
\begin{eqnarray}
\frac{d}{dv}L(g_v) & = & \frac{d}{dv} \int_0^{L_0} \beta(u,v) du
\nonumber \\
& = & \int_0^{L_0} \dr_v\beta du \nonumber \\
& \leq & \int_0^{L_0} (\tau_0+2\tau_1) \alpha(u,v) \beta(u,v) du 
\nonumber \\
& \leq & (\tau_0+2\tau_1) \left(\sup_{u\in [0, L_0]} \alpha(u,v)\right)
\int_0^{L_0} \beta(u,v) du~, \nonumber 
\end{eqnarray}
so:
\beq \label{var-Lg} 
\frac{d}{dv}L(g_v) \leq (\tau_0+2\tau_1) \exp((\tau_0+2\tau_1) L(g_v))
L(g_v) 
\eeq 
Now integrate this equation to obtain the required upper bound on $L_1$;
the upper bound on $L'_1$ is obtained in the same way, exchanging $u$
and $v$. Finally, the upper bound on the area comes from the upper
bounds on $L(g_v)$ and on $\sup_{u\in
[0, L_0]}\alpha(u,v)$ which we have found.
\epv

\bcr \label{cr-cont}
Let $x, y, z\in \Si$ be such that there exists an integral curve of $U$
(or of $-U$) of length at most $L_M$ going from  $x$ to
$y$, and an integral curve of $V$ (or of $-V$) of length at most 
$L_M$ going from $x$ to $z$. Then the integral curve of $U$ through
$z$ meets the integral curve of $V$ through $y$.
\ecr

\bpv
The intersections between those integral curves remain at bounded
distance as long as the lengths of the integral curves of $U$ and $V$
going from $x$ 
to $y$ and to $z$ remain below $L_M$, and the integral curves of $U$ and
of $V$ do not meet $\dr_{\III}\Si$ because $(\Si, I)$ is complete
(cf. section 2). 
\epv

\bpn{of lemma \ref{propag}}
Choose $u\in [0, L_g]$ and $t\in [-2T_0, 2T_0]$.
From corollary \ref{cr-cont}, if $T_0$ and $L_g$ are below a fixed
constant, then the integral curve of $U$ starting at $h_0(t)$ meets the
integral curve of $V$ starting at $g(u)$. Let $g_t(u)$ be their
intersection, and $A(u,t)$ be the area of the domain in $\Si$
bounded by $g$, $g_{[0,t]}(L_{g})$, $g_t$, $g_{[0,t]}(0)$. By
proposition \ref{controles}, there exists some $T_1>0$ such that, if
$L_g\leq T_1$ and $2T_0\leq T_1$, then, for all
$u\in [0, L_g]$:
$$ \left\|\frac{\dr g_t(u)}{\dr u}\right\|\leq 2~, ~~ 
\left\|\frac{\dr g_t(u)}{\dr t}\right\|\leq 2~. $$

By definition of $\delta_{g}$, there exists $u_0, u_1\in [0,L_g]$ such that
$\theta(g(u_0))-\theta(g(u_1))\geq \pi - \delta_{g}$. Then
$\theta(g(u_0))\geq \pi-\delta_g$ and $\theta(g(u_1))\leq
\delta_g$. To simplify the notations a little, we suppose that
$u_0<u_1$. 
For each $t\in [-2T_0, 2T_0]$, let:
$$ \theta_t:=\max(\pi-\theta(g_t(u_0)), \theta(g_t(u_1)))~. $$
We now suppose that $t\in [-2T_0, 2T_0]$ is such that, for all $s\in [0,
t]$, $\theta_s\leq 2C_1\eps$ and $\sigma_{g_s}\leq 2C_1\eps$, for some
constant $C_1$ on which more details will be given later. We will show
that, if $T_0$ is small enough, then this implies that $\theta_t\leq
C_1\eps$ and $\sigma_{g_t}\leq C_1\eps$, so that the same bounds apply
for all $t\in [-2T_0, 2T_0]$.

Let $C_t$ be the closed curve $g_0([u_0, u_1])\cup
g_{[0,t]}(u_1)\cup g_t([u_0, u_1])\cup
g_{[0,t]}(u_0)$. Let $W$ be the vector field on $C_t$ equal to $V$
over $g_0([u_0, u_1])\cup g_t([u_0, u_1])$ and to $U$ over
$g_{[0,t]}(u_1)\cup g_{[0,t]}(u_0)$. According to the Gauss-Bonnet
theorem, the total rotation of $W$ on $C_t$ (i.e. the integral of
$\langle \Nat W, JW\rangle$ plus the terms corresponding to the vertices) is
$-K_t$, where $K_t$ 
is the integral of $\Kt$ on the interior of $C_t$. But :
\begin{itemize}

\item The terms corresponding to $g_0([u_0, u_1])$ and to $g_t([u_0,
u_1])$ are bounded because of (\ref{UV2*}):
$$ \int_{g_0([u_0,u_1])}\|\Nat_UV\| \leq \int_{u_0}^{u_1} \tau_1\sin
\theta(g_0(u)) du \leq \tau_1 \sigma_{g_0} \leq 2\tau_1C_1\eps~, $$
$$ \int_{g_t([u_0,u_1])}\|\Nat_UV\| \leq \tau_1 \sigma_{g_t} \leq
2\tau_1C_1\eps~. $$

\item the terms corresponding to $g_{[0,t]}(u_1)$ and $g_{[0,t]}(u_0)$ are
bounded because of (\ref{VU2*}):
\begin{eqnarray}
\int_{g_{[0,t]}(u_0)}\|\Nat_VU\| & \leq & \int_0^t \tau_1 \left\| \frac{\dr
g_s(u_0)}{\dr s}\right\| \sin \theta(g_s(u_0)) ds \nonumber \\
& \leq & \int_0^t 2\tau_1 \theta_s ds \nonumber \\
& \leq & 4\tau_1 C_1\eps t~. \nonumber 
\end{eqnarray}

\item $K_t$ is bounded by:
\begin{eqnarray}
K_t & \leq & K_5 \int_{s=0}^t \int_{u=u_0}^{u_1} \left\| \frac{\dr
g_s(u)}{\dr s}\right\|  \left\| \frac{\dr
g_s(u)}{\dr u}\right\| \sin \theta(g_s(u)) du ds \nonumber \\
& \leq &  K_5 \int_{s=0}^t \int_{u=u_0}^{u_1} 4 \sin \theta(g_s(u)) du
ds \nonumber \\
& \leq &  4 K_5 \int_{s=0}^t \sigma_{g_s} ds~, \nonumber 
\end{eqnarray}
so that:
\beq \label{mKt}
K_t\leq 8 K_5 C_1 \eps t~. 
\eeq
\end{itemize}
Therefore:
$$ |\theta(g_0(u_0))-\theta(g_0(u_1)) +
\theta(g_t(u_1))-\theta(g_t(u_0))|\leq 4\tau_1 C_1\eps + 8 \tau_1
C_1\eps t + 8 K_5 C_1\eps t~. $$
Thus:
$$ |(\pi-\theta(g_t(u_0))) + \theta(g_t(u_1))| \leq
|(\pi-\theta(g_0(u_0))) + \theta(g_0(u_1))| + 4\tau_1 C_1 \eps (1+2t) +
8 K_5 C_1 \eps t~, $$
so that:
$$ |(\pi-\theta(g_t(u_0))) + \theta(g_t(u_1))| \leq \delta_g + 4\tau_1
C_1 \eps (1+2t) + 8 K_5 C_1 \eps t~. $$
This already shows that, if $T_0$ is such that:
$$ 1 + 4\tau_1 C_1 (1+2T_0) + 8K_5 C_1 T_0 \leq C_1~, $$
then $\theta_t\leq C_1\eps$. It remains to show that $\sigma_{g_t}\leq
C_1\eps$.

Equation (\ref{VU2*}) shows that:
\begin{eqnarray}
|\angle(U(g_t(u_0)), \Pi(g_{[0,t]}(u_0); U(g_0(u_0))))| & \leq & \tau_1
\int_{g_{[0,t]}(u_0)} \sin \theta(g_s(u_0)) ds \nonumber \\
& \leq & 2\tau_1 \int_0^t \theta_s ds \nonumber \\
& \leq & 2\tau_1 C_1 \eps t~, \nonumber
\end{eqnarray}
and, since $\theta_t\leq C_1\eps$ and $\theta_0\leq \eps$:
$$ |\angle(V(g_t(u_0)), \Pi(g_{[0,t]}(u_0); V(g_0(u_0))))| \leq 2\tau_1
C_1 \eps t + (C_1+1)\eps~. $$
On the other hand, (\ref{UV2*}) shows that, for all $u\in [0, L_g]$:
$$ |\angle(V(g_t(u)), \Pi(g_t([u_0, u]); V(g_t(u_0))))| \leq \tau_1
\sigma_t \leq 2\tau_1 C_1 \eps~, $$
while, for the same reason:
$$ |\angle(V(g_0(u)), \Pi(g_0([u_0, u]); V(g_0(u_0))))| \leq 2\tau_1 C_1
\eps~. $$ 
Finally, the same argument as in the proof of (\ref{mKt}) above shows
that the integral of $\Kt$ on the domain bounded by $g_0([u_0,u])$,
$g_{[0,t]}(u)$, $g_t([u_0,u])$ and $g_{[0,t]}(u_0)$ is at most
$8K_5C_1\eps t$. The Gauss-Bonnet theorem, applied to this domain,
therefore indicates that:
$$ |\angle(V(g_t(u)), \Pi(g_{[0,t]}(u); V(g_0(u))))| \leq 2\tau_1
C_1 \eps t + (C_1+1)\eps + 4\tau_1 C_1 \eps + 8K_5 C_1 \eps t~, $$
so that:
\beq \label{coV}
|\angle(V(g_t(u)), \Pi(g_{[0,t]}(u); V(g_0(u))))| \leq 8 C_1 \eps t
(\tau_1+K_5) + (C_1+1) \eps~.  
\eeq

Moreover, by (\ref{VU2*}):
\begin{eqnarray}
|\angle(U(g_t(u)), \Pi(g_{[0,t]}(u); U(g_0(u))))| & \leq &
\int_{g_{[0,t]}(u)} \|\Nat_VU\| \nonumber \\
& \leq & \int_0^t 2\tau_1 \sin \theta(g_s(u)) ds~, \nonumber
\end{eqnarray}
so that:
\begin{eqnarray}
\int_0^{L_g} |\angle(U(g_t(u)), \Pi(g_{[0,t]}(u); U(g_0(u))))| du  &
\leq & 2\int_0^{L_g} \int_0^t 2\tau_1 \sin \theta(g_s(u)) ds du
\nonumber \\ 
& \leq & 4\tau_1 \int_0^t \sigma_{g_s} ds \nonumber \\
& \leq & 8 C_1 \tau_1 \eps t~. \nonumber
\end{eqnarray}
But the definition of $\sigma_{g_t}$ shows that:
$$ \sigma_{g_t} \leq \int_0^{L_g} \sin \theta(g_0(u)) +
|\angle(U(g_t(u)), \Pi(g_{[0,t]}(u); U(g_0(u))))| + \coupeq +
|\angle(V(g_t(u)), \Pi(g_{[0,t]}(u); V(g_0(u))))| du~, $$
so by (\ref{coV}):
$$ \sigma_{g_t} \leq \sigma_{g_0} + 16 C_1 \tau_1 \eps t + 2 L_g 
(8 C_1 t (\tau_1+K_5) + C_1 + 1)\eps~. $$
But $|t|\leq 2T_0$, $L_g\leq T_0$ and $\sigma_{g_0}\leq 2 C_1\eps$, so
it is clear that there exists $C_1$ such that, if $T_0$ is small enough,
$\sigma_{g_t}\leq C_1\eps$. 

Using (\ref{coV}) once more then proves the lemma.
\epn

\section{Concave points}

We now turn to the proof of lemma \ref{dom}, which we recall here for
the reader's convenience. 

\bln{\ref{dom}}
Under the hypothesis of theorem \ref{thm-a}, $\dr_{\III}\Si$ has no
concave point. 
\eln

When $\phi:[-d,d]\times [0,d]\rightarrow \Si$ is a $(k,C)$-concave map,
and when $d'\leq d$, we call $\phi_{d'}$ the restriction of $\phi$ to
$[-d',d']\times [0, d']$; it is again a $(k,C)$-concave map. We also
call: 
$$ \dr_R \phi := \phi(\dr ([-d,d]\times [0,d])\cap \R\times \R_+^*)~. $$

Lemma \ref{dom} is a consequence of the following simpler lemma, whose
proof will be given below.

\blm \label{lm-normal}
Let $\phi:[-d,d]\times [0,d]$ be a concave map. There exist $\eps\in
(0,d]$ and $L>0$ such that, for any $\eps'\leq \eps$, there exists a
piecewise smooth curve $\ga$ of length at most $L$, which is an integral
curve of $U$ or $V$ on each interval where it is smooth, and which
goes from a point of $\im(\phi_{\eps'})$ to $\dr_R\phi_\eps$.
\elm

\bpn{of lemma \ref{dom}}
Suppose $x$ is a concave point of $\dr_{\III}\Si$. By definition, there
exists a concave map $\phi:[-d,d]\times [0,d]\rightarrow \Sib$ at $x$. 

Therefore, by lemma \ref{lm-normal}, there exists $\eps\in (0,d]$ such
that, for each $\eps'\leq \eps$, $\im(\phi_{\eps'})$ can be 
connected to $\dr_R\phi_\eps$ by a piecewise
smooth curve $\ga$ of length at most $L$ (for $\III$), which is an
integral curve of $U$ or $V$ on each smooth segment. 

By remark \ref{UV-I}, the length of $\ga$ for $I$ is at most
$L/\sqrt{K_2-K_1}$, so $\im(\phi_{\eps'})$ is at  
distance at most $L/\sqrt{K_2-K_1}$ from $\dr_R\phi_\eps$ for $I$. This
contradicts the fact that $(\Si, I)$ is complete.  
\epn

Now for the proof of lemma \ref{lm-normal}.
The proof is by contradiction, so we suppose that $\dr_{\III}\Si$ is
concave at a point $x_0$, with a $(k, C)$-concave map $\phi$ at
$x_0$, and such that there is no sequence of piecewise asymptotic curves
of bounded lengths starting from $\dr_R\phi_\eps$ (for some fixed
$\eps>0$) and ending arbitrarily close to $x_0$. We first state a remark
which will be used later on.

\bprop \label{rkW}
For $C_0, k_0>0$, there exists $C_1>0$ (depending on $C_0, k_0, \tau_0$ and
$K_5$) such that, if $\phi$ is a $(k_0, C_0)$-concave map, then the
``vertical'' curves $\phi(\{ x\}\times (0,d))$ have geodesic curvature
bounded by $C_1$. 
\eprop

\bpv
Let $X,Y$ be the vector fields and $v,l$ be the functions on $\im(\phi)$
such that:
$$ \dr_1\phi = v X ~, ~~ \dr_2\phi = l Y~. $$
Call $k:=\III(\Nat_YX,Y)$ and $\kappa:=\III(\Nat_XX,Y)$. Then:
$$ X.k = \III(\Nat_X\Nat_YX,Y) = \III(\Nat_Y\Nat_XX,Y) +
\III(R^{\Nat}_{X,Y}X, Y) + \III(\Nat_{[X,Y]}X,Y)~. $$
But $[vX,lY]=0$, so that:
$$ [X,Y] + \frac{dl(X)}{l} Y - \frac{dv(Y)}{v} X = 0~, $$
and therefore:
\beq \label{Xk}
X.k = Y.\III(\Nat_XX,Y) + \Kt + \frac{dv(Y)}{v} \III(\Nat_XX,Y) -
\frac{dl(X)}{l} \III(\Nat_YX,Y)~.
\eeq

Now, by definition of $\tau$:
$$ \Nat_{vX} (lY) - \Nat_{lY}(vX) - [vX,lY] = \tau(vX,lY)~, $$
so that:
$$ \Nat_XY - \Nat_YX + \frac{dl(X)}{l}Y - \frac{dv(Y)}{v}X = \tau~, $$
and therefore:
\beq \label{elv}
\frac{dl(X)}{l} = k + \III(\tau,Y)~, ~~\frac{dv(Y)}{v} = -\kappa -
\III(\tau, X)~.
\eeq
Using this and (\ref{Xk}) shows that:
$$ X.k = Y.\kappa + \Kt - \kappa (\kappa+\III(\tau,X)) -
k(k+\III(\tau,Y))~. $$

Now the definition of a convex map and the bounds on $\Kt$ and $\tau$
show that:
$$ |X.k|\leq C_0 + K_5 + C_0k_0(C_0k_0+\tau_0)+|k|(|k|+\tau_0)~. $$
This means that, if $k$ is large at a point $m$, then it remains large
in a neighborhood of $m$ in the integral curve of $X$ through
$m$. Equation (\ref{elv}) shows that $l$ would then vary a lot on this
curve, and this would contradict the definition of a convex map (because
$1\leq \|\dr_2\phi\|\leq C_0$).
\epv

The next point in the proof of lemma \ref{lm-normal} is to prohibit the
existence of asymptotic curves with $\delta$ small, using lemma
\ref{propag} and the following result: 

\bprop \label{cvx-qg}
For each $C>0$, there exist $\eps(C)>0$ such 
that, for any $\eps\leq \eps(C)$, if $\diam(\phi)\geq \eps$ and if
$c:[-4\eps,4\eps]\rightarrow \Si$ is a $Cy(c(0))$-quasi-geodesic, then
$c$ meets $\dr_R\phi_{\eps}$.   
\eprop

\bpv
Let $c:[0, L)\rightarrow \im(\phi)$ be an $\eps$-quasi-geodesic,
with either $L=\infty$ or $L\in \R_+^*$ and $\lim_L c \in \dr \phi$. 
Note $W$ the vector field on $\im(\phi)$ defined by:
$$ W:=\frac{\dr_1\phi}{\|\dr_1\phi\|}~. $$
Let $\alpha(t)$ be the angle between $W$ and $c'(t)$, and
$\alpha_0(t)$ the angle between $W$ and the parallel transport at
$c(t)$ of $c'(0)$ along $c([0,t])$.

By proposition \ref{rkW}:
$$ \|\Nat_{\dr_2\phi} W\|\leq C_1\|\dr_2\phi \|~, $$
while the definition of a $(k,C)$-concave map indicates that:
$$ k \leq \langle \Nat_WW, JW\rangle \leq Ck~. $$
As a consequence:
$$ k\cos \alpha(t) - C_1 |\sin \alpha(t)| \leq \langle \Nat_{c'(t)}W,
JW\rangle \leq Ck \cos \alpha(t)  + C_1|\sin \alpha(t) |~, $$
so that, by definition of a quasi-geodesic, if $\eps_0:=Cy(c(0))$, then:
\beq \label{valp0} 
k\cos \alpha_0(t) - C_1|\sin \alpha_0(t)| - (Ck+C_1)\eps_0 \leq
\coupeq \leq
\alpha_0'(t) \leq Ck\cos \alpha_0(t) +
C_1|\sin \alpha_0(t)| + (Ck+C_1)\eps_0~.
\eeq

We now suppose (without loss of generality) that $\alpha_0(0)\in [0,
\pi/2]$. 
Let $\alpha_1>0$ be the smallest positive number such that:
$$ k\cos(\alpha_1) - C_1 \sin(\alpha_1) - (Ck+C_1)\eps_0 \geq 0~. $$
$\alpha_1$ exists if $\eps_0$ is small enough (which happens if $\eps$
is small enough).
Equation (\ref{valp0}) indicates that, if $\alpha_0\in [\alpha_1,
\pi-\alpha_1]$ at a time $t$, then it remains there until the time $L$ where
$c$ leaves
$\im(\phi)$; moreover, $\alpha_0$ reaches $[\alpha_1,
\pi-\alpha_1]$ before a
fixed time $t_0$ (depending on $C, K$, etc) and, in the interval $[0,
t_0]$, it remains above $-c_0\eps_0$, where $c_0>0$ is a constant depending
also on $C$ and $k$.

Now it is easy to check that:
$$ \frac{\cos \alpha(t)}{C} \leq x'(t) \leq \cos \alpha(t)~, $$
$$ \frac{\sin \alpha(t)}{C} \leq y'(t) \leq \sin \alpha(t)~, $$
so that:
\beq \label{vx}
\frac{\cos \alpha_0(t)}{C} - \eps_0 \leq x'(t) \leq \cos \alpha_0(t) +
\eps_0~,
\eeq
\beq \label{vy}
\frac{\sin \alpha_0(t)}{C} - \eps_0 \leq y'(t) \leq \sin \alpha_0(t) +
\eps_0 ~. 
\eeq

As a
consequence of the lower bounds on $\alpha_0$ and on $y'$: 
$$ \forall t\in [0, L], ~ y(t)\geq y(0) - (c_0 + 1) t_0 \eps_0~, $$
so that $c$ must intersect $\dr_R\phi_\eps$ if $\eps$ is small enough. 
\epv

The same ideas also lead to the following statement, where we suppose
that $c'(0)$ is not too horizontal, instead of supposing that $y(c(0))$
is not too small. The proof is left to the reader.

\bprop \label{qg-angle}
There exists $C_{\angle}>0$ and $\eps_{\angle}>0$ such that, for any
$\eps\leq \eps_{\angle}$, if $\diam(\phi)\geq \eps$ and if
$c:[-4\eps,4\eps]\rightarrow \Si$ is a $\alpha_0$-quasi-geodesic with
$C_{\angle}\alpha_0\leq \angle(\dr_1\phi, c'(0)) \leq \pi/2$, then
$c$ meets $\dr_R\phi_{\eps}$.   
\eprop

The situation is simpler if $c$ is a geodesic instead of a
quasi-geodesic:

\bprop \label{c-geod}
There exist $C_3, \eps_3>0$ as follows.
Let $\eps\leq \eps_3$, and suppose that $\diam(\phi)\geq
\eps$. 
Let $g:[a,b]\rightarrow \im(\phi_{2\eps})$ be a maximal geodesic segment
in $\im(\phi_{2\eps})$, with $a<0<b$ and $g(0)\in \im(\phi_\eps)$. Then
either $g(a)$ or $g(b)$ is on $\dr_R\phi_{2\eps}$, and $L(g)\in
[\eps/C_3, C_3\eps]$. If $g'(0)$ is parallel to $\dr_1\phi$, then
both $g(a)$ and $g(b)$ are on $\dr_R\phi_{2\eps}$.
\eprop

\bpv
The proof is similar to that of proposition \ref{cvx-qg}; we call
$\alpha_0(t):=\angle(W, g'(t))$, and we suppose that $\alpha_0(0)\in [0,
\pi/2]$. Then:
$$
k\cos \alpha_0(t) - C_1|\sin \alpha_0(t)| \leq
\alpha_0'(t) \leq  Ck\cos \alpha_0(t) +
C_1|\sin \alpha_0(t)|~.
$$
It is the clear that $\alpha_0$ will soon become positive;
moreover, if we let $x(t):=x(g(t))$ and $y(t):=y(g(t))$, then:
$$
\frac{\cos \alpha_0(t)}{C}\leq x'(t) \leq \cos \alpha_0(t)~,
$$
$$
\frac{\sin \alpha_0(t)}{C}\leq y'(t) \leq \sin \alpha_0(t)~. 
$$
The second equation indicates that $y(t)$ remains positive while
$g(t)\in \im(\phi_{\eps})$, and both equations taken together again
show that $g$ intersects $\dr_R\phi_{\eps}$ after time at most
$C_3'\eps$ for some $C_3'>0$. 

The same equations apply for the segment of $g$ where $t\leq 0$; after a
bounded time, either $g$ will have intersected $\dr\phi_{2\eps} \setminus
\dr_R\phi_{2\eps}$, or $y'$ will vanish. The same argument as above then
shows that, in both cases, $t\mapsto g(-t)$ will meet $\dr\phi_{2\eps}$
after a time at most $C_3''\eps$ for some $C_3''>0$. This proves the
upper bound on the length of $g$.
The corresponding lower bound comes from the distance between
$\im(\phi_{\eps})$ and the part of $\dr_R\phi_{2\eps}$ that $g$ can
intersect for $t>0$.

Finally, the case where $g'(0)$ is parallel to $\dr_1\phi$ is obtained
by applying twice the argument for $t>0$, which can be used in this case
also for $t<0$ because $-g'(0)$ is also directed towards the increasing
values of $y$.
\epv

From now on, we consider an integral curve $g:I\rightarrow \Si$ of
$U$, where $I$ is an interval, either of the form $[0, t_M]$ or
$\R_+$. $\eps$ is a fixed positive number, on which more details are
given below.

\bdf 
Let $t\in I$; call $\ga_t$ the maximal geodesic segment directed by
$V(g(t))$. Call $E_g$ the subset of $I$ containing all $t$ such that $\ga_t$
intersects $\dr_R\phi_\eps$ on both sides at finite 
distance. For $t\in E_g$, call $\Om_t$ the connected component of
$\im(\phi_\eps)\setminus \ga_t$ which does not contain
$x_0$ in its boundary. 
\edf

\bprop \label{dom-dec}
If $t\in E_g$ and $g'(t)$ is towards the interior of $\Om_t$, then, for
all $t'\in I$ with $t'\geq t$, $t'\in E_g$, and $\Om_{t'}\subset \Om_t$. 
\eprop

\bpv
Note that, if $\eps$ is small enough, then, for any $t'\in E$, if
$t''>t'$ is close enough to $t'$, then $\ga_{t'}\cap
\ga_{t''}=\emptyset$. This comes from (\ref{UV2*}) and from corollary
\ref{cr-ido}.  This immediately implies that $(\Om_t)$ is a decreasing
family of subsets of $\im(\phi_{\eps})$.
\epv

\vspace{0.3cm}
\centerline{\psfig{figure=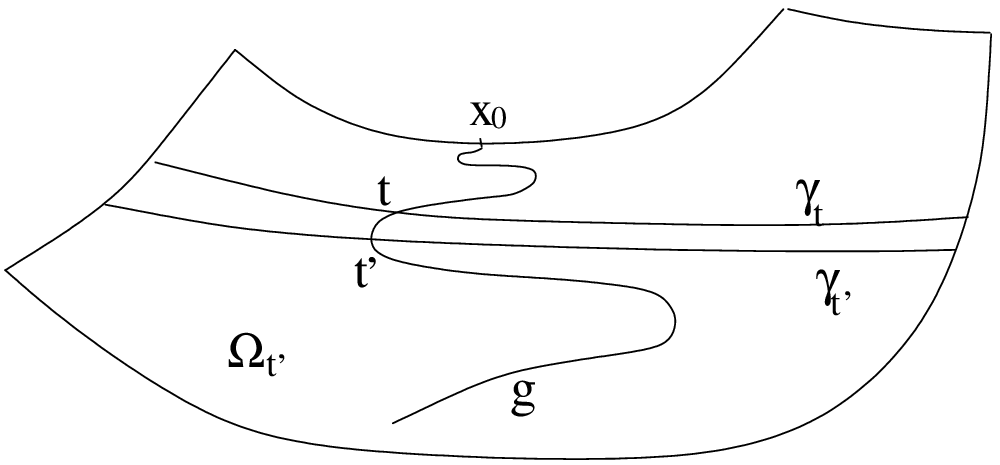,height=4cm}}
\centerline{\bf Figure 5.1} \vspace{0.3cm}

This is now used to prove that there exists an asymptotic curve going
from $x_0$ to $\dr_R\phi_\eps$.

\bprop \label{ex-ciu} 
If $\eps$ is small enough, there exists an integral curve $g:\R_+\rightarrow
\im(\phi_\eps)$ of $U$ (or of $V$) such that $g(0)\in \dr_R\phi_\eps$
and that $\lim_{t\rightarrow \infty} g(t)=x_0$. 
\eprop

\bpv
Fix $\eps'>0$. For $z\in (0, \eps')$, let $c_z$ be the maximal integral
curve of $V$ in $\im(\phi_{\eps'})$ containing $\phi(0,z)$. 
We consider two cases.

{\bf 1.} There exists $\eps>0$ and a sequence $z_n\rightarrow 0$ such
that, for each $n$, $c_{z_n}$ has one end on $\dr_R\phi_\eps$. 

If there exists $n$ such that $\overline{c_{z_n}}\ni x_0$, then the
proposition is proved. Otherwise, call $D_n$ the connected component of
$\im(\phi_{\eps'})$ which does not contain $x_0$ in its closure.
Let $D:=\cup_n D_n$.
Since the $c_{z_n}$ are integral curves of $V$, they are pairwise
disjoint (except when they coincide), so (maybe after taking a
subsequence of $(z_n)$) $(D_n)$ is an increasing sequence. Since
$c_{z_n}\ni \phi(0,z_n)\rightarrow x_0$, it is then
not difficult to prove that $\dr D$ contains an integral curve of $V$
connecting $\dr_R\phi_{\eps'}$ to $x_0$.

{\bf 2.} For all $\alpha>0$, there exists $z_\alpha>0$ such that, for
$z\leq z_\alpha$, $c_z$ remains in $\im(\phi_\alpha)$ and has both ends
on $\dr\phi_\alpha\setminus \dr_R\phi_\alpha$. 

Call $m_z$ a point of
$c_z$ where $y$ is maximal. Let $g_z:[0, L_z]\rightarrow
\im(\phi_{\eps'})$ be the maximal integral curve of $U$ (or $-U$) with
$g_z(0)=m_z$ and $g'(0)$ directed towards the increasing values of
$y$. By definition of $m_z$, $V(g(0))$ is parallel to
$\dr_1\phi$. Therefore, by proposition \ref{c-geod}, the geodesic directed
by $V(g(0))$ meets $\dr_R \phi_{\eps'}$ on both sides. With the
notations above, this indicates that $0\in E_{g_z}$, so that, by
proposition \ref{dom-dec}, $E_{g_z}=[0, L_z]$. Thus $g_z(L_z)\in
\dr_R\phi_{\eps'}$. 

The proof then proceeds as in case (1.) above, because the $(g_z)$ are
disjoint and, after taking a sequence $(z_n)\rightarrow 0$, they
converge to an integral
curve of $U$ connecting $\dr_R\phi_{\eps'}$ to $x_0$.
\epv

Moreover, the rate of decrease of the area of $\Om_t$ is bounded by
$\sin \theta(g(t))$:

\bprop \label{min-aire}
There exists $\lambda_4>0$ such that, if $t\in E$ is such that $g(t)\in
\im(\phi_{\lambda_4\eps})$, then:
$$ \left| \frac{\dr \area(\Om_t)}{\dr t} \right| \geq \lambda_4 \sin
\theta(g(t)) \eps~, $$
and:
$$ \liminf_{t'\rightarrow t^+} \frac{d(\ga_{t'}, \ga_t)}{t'-t} \geq
\lambda_4 \sin \theta(g(t))~. $$ 
\eprop

\bpv
This is again a consequence of corollary \ref{cr-ido}, along with
(\ref{UV2*}), which bounds the rate of variation of the direction of $V$
along $\ga$.
\epv

For each $t\in \R_+$, we let $h_t:[-T_0,T_0]\rightarrow \Si$ be the
integral curve of $V$ with $h_t(0)=g(t)$. A direct consequence of the
previous proposition is that the integral of $\sin\theta$ on $g$ is
finite, and this will be used now to show that many $h_t$ are
quasi-geodesics. 

\bprop \label{qg-x0}
If $\eps$ is small enough, there exists $C>0$ such that, for each
$t_0\geq C$, there exists $t\geq t_0$ such that $h_t$ intersects
$\dr_R\phi_\eps$. 
\eprop

\bpv
Since the integral of $\sin \theta$ on $g$ is finite, equation
(\ref{UV2*}) shows that:
$$ \lim_{t,t'\rightarrow \infty} \angle(V(g(t')), \Pi(g_{|[t,t']};
V(g(t)))) = 0~. $$
One can therefore define a parallel vector field on $V_0$ as:
$$ V_0(g(t)) := \lim_{t'\rightarrow \infty} \Pi(g_{|[t,t']}; V(g(t')))~,
$$
and, by (\ref{UV2*}):
$$ |\angle(V_0(g(t)), V(g(t)))|\leq \tau_1 \int_t^\infty \sin \theta(g(s))
da \stackrel{t\rightarrow \infty}{\longrightarrow} 0~. $$

The same works for $W$ because $\lim_{\infty}g=x_0$; set:
$$ W_0(g(t)) := \lim_{t\rightarrow \infty} \Pi(g_{|[t,t']}; W(g(t')))~,
$$
and then:
$$ \lim_{t\rightarrow \infty} \angle(W_0(g(t)), W(g(t))) = 0~. $$

Let $\alpha_0:=\angle(W_0, V_0)$; we suppose (without loss of
generality) that $\alpha_0\in [0, \pi/2]$. The proof will proceed
differently according to whether $\alpha_0=0$ or $\alpha_0>0$.

If $\alpha_0>0$, remark that, since $\lim{\infty}g=x_0$, an elementary
argument (as e.g. in the proof of proposition \ref{cvx-qg}) shows that:
$$ \lim_{t\rightarrow \infty} \int_t^{t+T_0} \cos \theta(g(s)) ds = 0~. $$
Thus, for any fixed $t_0\in \R_+$ and $\lambda>0$, there exist $u,v\in
\R_+$ such that:
\begin{itemize}
\item $t_0\leq u\leq v\leq u+T_0$;
\item $\theta(g(u))\leq \lambda$ and $\theta(g(v))\geq \pi-\lambda$;
\item $\int_u^v\sin \theta(g(s)) ds \leq \lambda$;
\item $\angle(\dr_1\phi, V)\geq \alpha_0/2$ at $g(u)$.
\end{itemize}  
Then:
$$ \delta_{g_{|[u,v]}} \leq 2\lambda~, ~~\sigma_{g_{[u,v]}}\leq
\lambda~. $$
According to lemma \ref{propag}, $h_u$ is a quasi-geodesic; proposition
\ref{qg-angle} then indicates that, if $\lambda$ and $\eps$ are small
enough, $h_u$ intersects $\dr_R\phi_\eps$.

Consider now the case where $\alpha_0=0$. By proposition
\ref{min-aire}, there exists $c>0$ so that, for $t$ large enough:
$$ y(g(t))\geq c\int_t^\infty \sin\theta(g(s)) ds~. $$
But, again:
$$ \lim_{t\rightarrow \infty} \int_t^{t+T_0} \cos \theta(g(s)) ds = 0~,
$$
and the same argument as in the case $\alpha_0=0$ leads to the
conclusion, but with proposition \ref{qg-angle} replaced by proposition
\ref{cvx-qg} to shows that $h_u$, which is a quasi-geodesic, actually
intersects $\dr_R\phi_\eps$ if $\eps$ is small enough.
\epv

The proof of lemma \ref{lm-normal} obviously follows, because the
conclusion of the previous proposition contradicts the hypothesis, made
above, that the conclusion of lemma \ref{lm-normal} does not hold.

\section{The boundary is convex}

This section contains the proof that $\Sigma$,
with $\III$ and $\Nat$, is convex in the sense of definition
\ref{df-cvx}. 
 
\bln{\ref{ccv-cvx}}
If $\Si$ has no concave point, then it is convex.
\eln

We need to make
normal deformations of curves, while controling their curvature. 
Some of the tools
needed here will be used again in the next section, to prove that convex
surfaces have bounded area. 

From now on, whenever we consider a smooth, convex curve $g$, we suppose
that it is parametrized in such a way that $J_{\III}g'$ is oriented
towards the convex side of $g$.

\bprop \label{var-k}
Let $g:[0, L]\rightarrow \Si$ be a smooth curve with $\|g' \|\equiv
1$. Let $l:[0, L]\rightarrow \R$ be a smooth function. The
first order variation of the curvature $\kappa$ of $g$ in a 
normal deformation of $g$ which is defined at $g(s)$ by the vector $l(s)
Jg'(s)$ is:
$$ 
\kad = l (\Kt+\kappa
(\kappa+\tau(g'(s))))+X.X.l-X.(l\tau(Jg'(s)))~.
$$
\eprop

\bpv
By linearity, it is enough to prove this proposition when $l$ is
positive. 
Let $(g_{s})_{s\in [0,1]}$ be a one parameter family of curves such that
$g_0=g$ and that:
$$ \dr g_s(t)/\dr s\equiv l Jg'(s)/\| g' \|~. $$
To simplify somewhat the notations, we call $v_s$ the speed of $g_s$,
that is:
$$ v_{s}(t):=\left\|\dr g_{s}(t)\over \dr t\right\|_{\III}~. $$
We also call $X$ the unit vector along $g'_s(t):=\dr g_s(t)/\dr t$, 
and $Y:=JX$. Therefore:
$$
\dr_{s}(g_{s}(t)) = l_s Y
$$
for some function $l_s(t)$ such that $l_0=l$. Moreover,
$\kappa=\kappa_{s}(t)=\III(\Nat_{X}X,Y)$.

By definition of the torsion $\tau$ of $\Nat$,
\begin{equation}
\label{eq-torsion}
\Nat_{X}Y-\Nat_{Y}X-[X,Y]=\tau~,
\end{equation}
while, since $(s,t)$ define a coordinate system on a domain of $\Si$:
\begin{equation}
\label{eq-torsion2}
[vX,lY]=0~.
\end{equation}
Set $\lambda:=\ln(l)$ and $\nu:=\ln(v)$, the previous equation becomes:
\begin{equation}
\label{eq-torsion3}
[X,Y]+(X.\lambda)Y-(Y.\nu)X =0~.
\end{equation}
Let $\tau_X:=\langle \tau, X\rangle$ and $\tau_Y:=\langle \tau,
Y\rangle$. Then:
\begin{eqnarray}
Y.\kappa & = & \III(\Nat_{Y}\Nat_{X}X,Y) \nonumber \\
& = & \Kt + \III(\Nat_{X}\Nat_{Y}X-\Nat_{[X,Y]}X,Y) \nonumber  \\
& = & \Kt + \III(\Nat_X(\Nat_X Y -[X,Y] -\tau) + (X.\lambda)\Nat_Y X -
(Y.\nu)\Nat_X X, Y) \nonumber \\
& = & \Kt + \III(\Nat_X(-\kappa X + (X.\lambda)Y - (Y.\nu)X -\tau_X X
-\tau_Y Y), Y) + \nonumber \\
& & + (X.\lambda)\III(\Nat_X Y - [X,Y] -\tau, Y) - \kappa
(Y.\nu) \nonumber \\
& = &  \Kt - \kappa^2 + X.X.\lambda -\kappa (Y.\nu) - \kappa \tau_X
-X.\tau_Y + (X.\lambda)^2 - (X.\lambda)\tau_Y - \kappa (Y.\nu)~, \nonumber
\end{eqnarray}
and, since $Y.\nu=-\kappa-\tau_{X}$ by (\ref{eq-torsion}) and
(\ref{eq-torsion3}): 
\begin{eqnarray}
Y.\kappa & = & \Kt + X.X.\lambda + \kappa^2 + \kappa \tau_X - X.\tau_Y +
(X.\lambda)^2 - (X.\lambda)\tau_Y \nonumber \\
& = & \Kt + X.X.\lambda + (X.\lambda)^{2}+
\kappa (\kappa+\tau_X)-(X.\tau_Y)-\tau_Y(X.\lambda)~, \nonumber
\end{eqnarray}
so that:
\begin{eqnarray}
\dr_{s} \kappa & = & l (Y.\kappa)  \nonumber \\
& = & l (\Kt+\kappa (\kappa+\tau_X))+X.X.l-X.(l\tau_Y)~, \nonumber 
\end{eqnarray}
which is the formula we need.
\epv

As a consequence, we find an inequality:

\bcr \label{cr-df}
The rate of variation of $\kappa$ is bounded from below by:
$$ \kad \geq \frac{l}{4}((4K_4-\tau_0^2) + \tau(Jg')^2) + l'' -
(l\tau(Jg'))'~. $$ 
If the curvature of $g$ is bounded from above by $\kappa_M$, then $\kad$
is also bounded from above:
$$ \left[ \frac{l}{4}((4K_4-\tau_0^2) + \tau(Jg')^2) + l'' -
(l\tau(Jg'))' \right] + l \left((K_5-K_4) + \kappa_M^2 + \kappa_M\tau_0 +
\tau_0^2\right) \geq \kad~. $$ 
\ecr

\bpv
For any $\kappa\in [0, \kappa_M]$, we have:
$$
\kappa_M(\kappa_M+\tau_0) \geq \kappa (\kappa+\tau_{X}) \geq
\frac{\tau_{X}^{2}}{4}~. 
$$ 
Moreover, $\tau_{X}^{2}+\tau_{Y}^{2}=|\tau |^{2}\leq
\tau_{0}^{2}$ and $K_4\leq \Kt\leq K_5$, so that:
$$
K_5 + \kappa_M(\kappa_M+\tau_0) \geq \Kt+\kappa (\kappa+\tau_{X}) \geq
K_{4}-\frac{\tau_{0}^{2}}{4}+\frac{\tau_{Y}^{2}}{4}~, 
$$
and the corollary follows.
\epv

This means that trying to deform curves leads to a natural question on
solutions of differential equations; we will need the following
proposition.

\bprop \label{edo}
For each $\eps>0$ small enough, there exists $M_0\geq 1$ and $S_1\geq S_0>0$
such that, if $u\in C^\infty(\R, 
[-1/\eps, 1/\eps])$, there exists $s_0, s_1\in [S_0, S_1]$ and $y\in
C^\infty([0, s_1], [0, M_0])$ such that:
\beq \label{edo1} 
y'' = (yu)'-(\eps+\frac{u^{2}}{4}) y~,
\eeq
with:
$$ y(s)\in [1, M_0] ~ \mbox{for} ~ s\in [0, s_0]~, $$
$$ y(0)=1, \; y'(0)=u(0)+4~, $$
$$ y(s_0)=1, \; y'(s_0)\leq u(s_0)+4~, $$
$$ y(s_1)=0, \; y'(s_1)\leq 0~. $$
Moreover, for each $s\in [0, s_0]$, $|y'(s)|\leq M_0$.
\eprop

\bpv
Let $z=y'-yu$. The relation (\ref{edo1}) becomes:
\begin{equation}
\label{sys2b} \left\{ \begin{array}{ccl}
y' & = & yu + z
\\
z' & = & -(\eps+u^{2}/4) y~.
\end{array}
\right.
\end{equation}
Let:
\begin{equation}
\label{def3b}
X:=\left(
\begin{array}{c}
y
\\
z
\end{array}
\right) \; \mbox{and} \; 
m:=\left(
\begin{array}{cc}
u & 1
\\
-(\eps+u^{2}/4) & 0
\end{array}
\right)~.
\end{equation}
The relation now is:
\begin{equation}
\label{sys3b}
X'(s) = m(s) X(s)~,
\end{equation}
so, for $s\geq 0$:
\begin{equation}
\label{def4}
X(s) = \exp(s M(s))X(0) = \exp \left( s\left(
\begin{array}{cc}
F(s) & 1
\\
-(\eps+\Ft^{2}/4) & 0
\end{array}
\right) \right) X(0)~,
\end{equation}
with:
$$
F(s)=\frac{1}{s}\int_{0}^{s} u(r)dr \; , \; \Ft(s) = 
\sqrt{\frac{1}{s}\int_{0}^{s}u(r)^{2} dr}~.
$$ 

Now the eigenvalues of $M(s)$ are the roots of:
$$
X(X-F)+\left(\eps+\frac{\Ft(s)^{2}}{4}\right)=0~.
$$ 
From the Cauchy-Schwarz theorem, $\Ft(s)^2\geq F(s)^2$,
so $4\eps+\Ft^{2} - F^{2} \geq 4\eps > 0$, 
and the eigenvalues are $\mu_{\pm}=\alpha \pm i\beta$, 
where:
$$ |\alpha|=\frac{|F|}{2} \leq \frac{1}{2\eps}~, ~~ |\beta| =
\sqrt{\eps + \frac{\Ft^2-F^2}{4}} \geq \sqrt{\eps}~. $$
The associated eigenvectors are $v_{\pm}=(1,-\alpha\pm i\beta)$. 

Suppose now that $X(0)=(1,4)$. Then
$y'(0)>0$, and $y(s)$ remains above $1$ until after time $s_0>0$. But
$\beta \geq \sqrt{\eps}$, so it is clear that, after a time $s_1\leq 
\pi/\sqrt{\eps}$, $y(s)$ will become negative,
and this provides an upper bound for $s_0$ and for $s_1$: 
\begin{equation}
\label{maj-s0}
s_0\leq s_1\leq S_1 := \frac{\pi}{\sqrt{\eps}}
\end{equation}

The decomposition of $(1,4)$ on the basis
$(v_{+},v_{-})$ is: 
\begin{equation}
\label{dec-X0}
(1,4)=\alpha_{+}v_{+}+\alpha_{-}v_{-}=\left(\frac{1}{2} -
i\frac{\alpha+4}{2\beta}\right) v_{+} +
\left(\frac{1}{2}+i\frac{\alpha+4}{2\beta}\right) v_{-}~, 
\end{equation}
and this gives an upper bound for $\alpha_{\pm}$: 
\begin{equation}
\label{maj-mu}
|\alpha_{\pm}|^{2} \leq \frac{1}{4} + \frac{(|\alpha|+4)^2}{4\beta^2}
\leq \frac{1}{4} + \frac{1}{4\eps}\left( \frac{2}{\eps}+4\right)^2~. 
\end{equation}
Moreover, $\beta$ also has an upper bound because $|u|\leq 1/\eps$, and
this gives a lower bound $S_0$ for $s_0$. 

Using the upper bound on $s_0$ (equation
(\ref{maj-s0})) and on $\alpha$ ($|\alpha| \leq 1/2\eps$), we see that, for
any $s\in [0, t_0]$:
\begin{eqnarray}
y(s) & \leq & |X(s)| \nonumber \\
& \leq & (|\alpha_{+}|+|\alpha_{-}|) \exp(s_0) \nonumber \\\
& \leq & \sqrt{1 + \frac{1}{\eps}\left(\frac{2}{\eps}+4
\right)^2}
\exp(\frac{\pi}{\sqrt{\eps}}) \nonumber~.
\end{eqnarray}
In addition, $y(s)>0$, so that $z'(s)<0$ for $s\in [0, s_{0}]$,
and it follows that $z(s_{0}) \leq z(0) = 4$. The function $y$ therefore
verifies the conclusions of proposition \ref{edo}.
\epv

\bcr \label{cr-defo}
For any $\kappa_M>0$, there exist $L_1>0, C_1>0$ and $M_1\geq 1$ as
follows. Let $\kappa_m\in [0, \kappa_M]$, and let $g_0:[0, L]\rightarrow
\Si$ be a convex curve parametrized at speed one, with
curvature $\kappa\in [\kappa_m, \kappa_M]$, and with $L\in (0,
L_1]$. Then there 
exists $T_1>0$ and a deformation $(g_t)_{t\in [0, T_1]}$ such that:
\begin{enumerate}
\item for each $t\in [0, T_1]$, $g_t$ is a convex curve with curvature
$\kappa\in [\kappa_m+t/C_1, \kappa_M+t C_1]$; 
\item for any $t\in [0, T_1]$ and $s\in [0, L]$, $(\dr_t g_t)(s)$ is
parallel to the unit normal $N$ to $g_t$, and $\langle (\dr_t g_t)(s),
N\rangle \in [1, M_1]$;
\item for each $t\in [0, T_1]$ and $s\in [0, L]$, $|\dr\kappa/\dr t|\leq
C_1$. 
\end{enumerate}
\ecr

\bpv
Let $\eps_0:=4K_4-\tau_0^2$, so that $\eps_0>0$ by lemma \ref{Sit}. By
proposition \ref{edo}, if $L$ is small enough, there exists a function
$l:[0, L]\rightarrow [1, M_0]$ such that, on $[0, L]$:
$$ l'' + \left(\frac{\eps_0}{2} + \frac{\tau(Jg_0')^2}{4}\right) l - (l
\tau(Jg_0'))' = 0~, $$
so that:
$$ l'' + \left(\eps_0 + \frac{\tau(Jg_0')^2}{4}\right) l - (l \tau(Jg_0'))'
= \frac{\eps_0 l}{2} \in \left[\frac{\eps_0}{2}, \frac{\eps_0
M_0}{2}\right]~. $$ 
For $t$ small enough and $s\in [0, L]$, set:
$$ g_t(s) = \exp_{g_0(s)}{tN(s)}~, $$
where $N(s)$ is the unit normal to $g_0$ at $g_0(s)$ towards the convex
side of the complement. This defines, for $t$ small enough, a smooth
curve $g_t$. 

Let $R_0:=(K_5-K_4)+\kappa_M^2+\kappa_M\tau_0+\tau_0^2$. 
Corollary \ref{cr-df} shows that, for $t=0$:
$$ \kad\in \left[\frac{\eps_0}{2}, \frac{M_0(\eps_0+R_0)}2\right]~. $$ 
Then, by compactness, for $t$ small enough,
$\kad \in [\eps_0/4, M_0(\eps_0+R_0)]$, and the corollary follows.
\epv

\bcr \label{cr-dt}
There exists $\eps_2>0$, $C_2>0$ and $L_2>0$ as follows. Suppose that
$L\in (0, L_2]$, and let $g:[0, L]\rightarrow \Si$ be a geodesic
segment. Let $d_0:=d_{\III}(g([0, L]), \dr_{\III}\Si)$, and suppose that
$d_0\leq \eps_2$. Suppose moreover that $d_{\III}(g(0),
\dr_{\III}\Si)\geq C_2 d_0$ and that $d_{\III}(g(L), \dr_{\III}\Si)\geq
C_2 d_0$. Then $(\Si, \III, \Nat)$ has a concave point.
\ecr

\bpv 
Apply corollary \ref{cr-defo} recursively to obtain a $(k, C)$-concave
map $\phi:[-d,d]\times [0, d]\rightarrow \Sib$ such that
$\phi([-d,d]\times \{ d\})$ is a segment of $g$.
\epv

\bpn{of lemma \ref{ccv-cvx}}
Suppose that $(\Si, \III, \Nat)$ is not convex. Then, by definition
\ref{df-cvx}, there exists a sequence of geodesic segments $\ga_n:[0,
L]\rightarrow \Si$ such that $(\ga_n(s))_{n\in \N}$ converges in $\Sib$
for each $s\in [0, L]$, with $\lim_{n\rightarrow \infty} \ga_n(s')\in
\dr_{\III}\Si$ 
for some $s'\in (0, L)$ but $\lim_{n\rightarrow\infty} \ga_n(0)\in
\Si$. Set: 
$$ s_0 := \inf \{ s\in [0, L] \; | \; \lim_{n\rightarrow \infty}
\ga_n(s)\in \dr_{\III}\Si \}~. 
$$ 

Remark that there exists $\eps_3\in (0, \max(s_0/3, L_2/2))$such
that, for each $n\in \N$, $B(\ga_n(s_0), 3\eps_3)\setminus \ga_n([0,
L])$ has at least two connected components, one of which is a half-disk
which does not meet $\dr_{\III}\Si$. Otherwise, each ball centered at
$\ga_n(s_0)$ would meet $\dr_{\III}\Si$ on each side of $\ga_n$ for each $n\in
\N$, and then there could be no path joining $\ga_n(0)$ to $\ga_n(L)$
(for $n$ large enough) in $\Si$, a contradiction. We suppose that the
half-disk which does not meet $\dr_{\III}\Si$ is always on the same side of
$\ga_n$ as $J\ga'_n(s_0)$.

By definition of $s_0$, $\ga_n([0, s_0-\eps_3])$ remains in a compact
subset of $\Si$, so that there exists $\eps_4>0$ so that, for each $n\in
\N$, $d_{\III}(\ga_n([0, s_0-\eps_3]), \dr_{\III}\Si)\geq 2\eps_4$. We call:
$$ \Om := \{ x\in \Si \; | \; \exists n\in \N, d_{\III}(x, \ga_n([0,
s_0-\eps_3)))\leq \eps_4 \}~, $$ 
so that $d_{\III}(\Om, \dr_{\III}\Si)\geq \eps_4$.

Let $\theta\in (0, \pi)$. Call $s(\theta)$ the supremum of all $s\in [0,
s_0]$ such that, for any $n\in \N$ and any $t\in [0, s]$, the maximal
geodesic starting from $\ga_n(s)$ with $\angle(\ga'_n(s),
g'(0))=\theta-\pi$ does not reach $\dr_{\III}\Si$ before time at least
$\eps_3$. Then, clearly:
$$ \limsup_{\theta\rightarrow 0^+} s(\theta)=s_0~. $$ 
For $s\in [s_0-\eps_3, s_0]$, $s<s(\theta)$, and $n\in \N$, the geodesic
segment $\gb$ starting from $\ga_n(s)$ with $\angle(\ga'_n(s),
g'(0))=\theta$ also does not reach $\dr_{\III}\Si$ before time at least
$\eps_3$, because it remains in a half-disk bounded by $\ga_n$ and of
radius $3\eps_3$. Let $g_{n, \theta, s}:[-\eps_3, \eps_3]\rightarrow
\Si$ be the geodesic segment with $g_{n, \theta, s}(0)=\ga_n(s)$ and 
$\angle(\ga_n'(s), g_{n, \theta, s}'(0))=\theta$.

\vspace{0.3cm}
\centerline{\psfig{figure=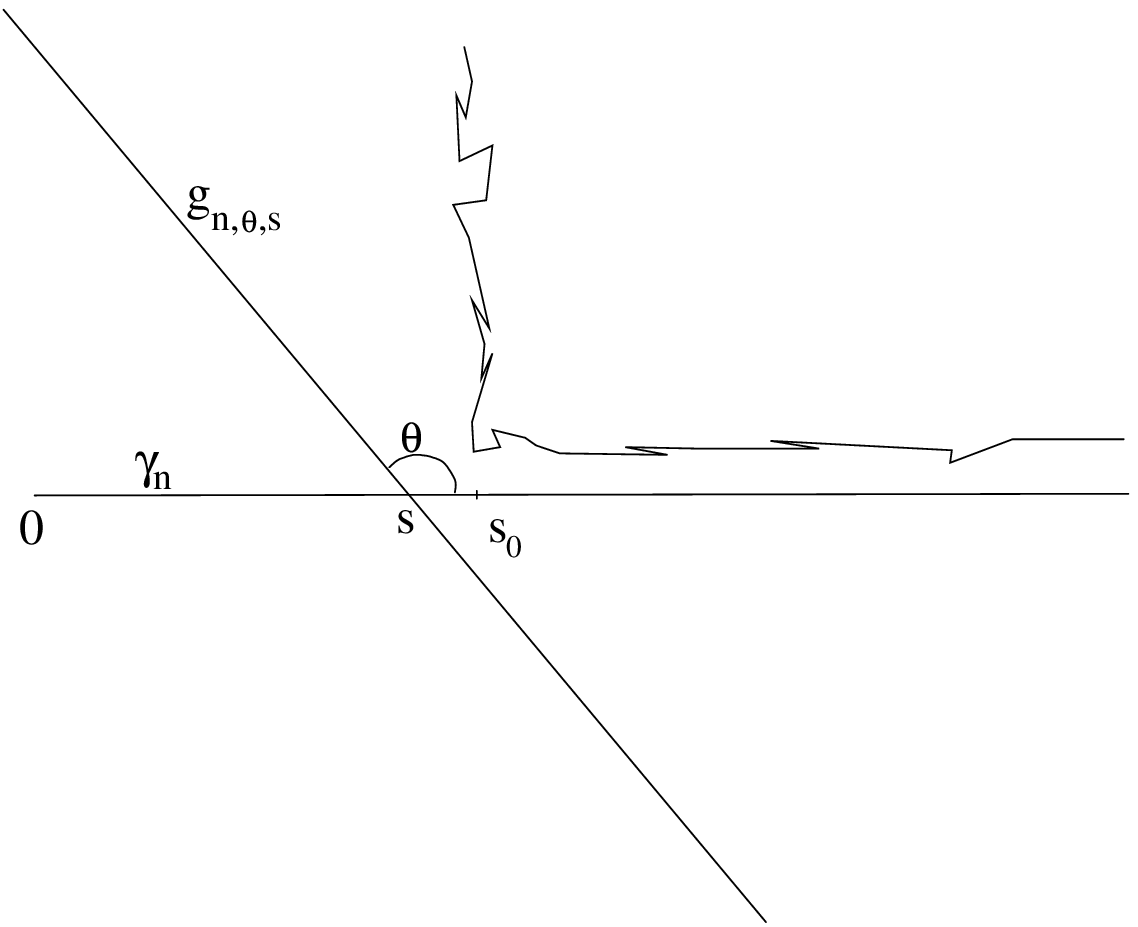,height=6cm}}
\centerline{\bf Figure 6.1} \vspace{0.3cm}

Now it is easy to check that, if $\theta$ is smaller than some fixed
$\theta_4$, then, for any $n\in \N$ and any $s\in [s_0-\eps_3,
s(\theta)]$, $d_{\III}(g_{n, \theta, s}(-\eps_3), \ga_n([0,
s_0 - \eps_3]))\leq \eps_4$; then $g_{n, \theta, s}(-\eps_3)\in \Om$, so
that $d_{\III}(g_{n, \theta, s}(-\eps_3), \dr_{\III}\Si)\geq \eps_4$.

On the other hand, one can check that $d_{\III}(g_{n, \theta_4, s}(\eps_3),
\ga_n)$ is bounded below by some fixed $\eps_5>0$ depending only on
$\theta_4$ (and $K_5, \tau_0$). Since $g_{n, \theta_4, s}([0, \eps_3])$
remains in a half-disk bounded by $\ga_n$ and of radius $3\eps_3$, this
shows that $d_{\III}(g_{n, \theta_4, s}(\eps_3), \dr_{\III}\Si)\geq
\eps_5$. Finally, by definition of $s(\theta_4)$, we can choose $n$ and
$s$ so that $d_{\III}(g_{n, \theta_4, s}([-\eps_3,\eps_3]),
\dr_{\III}\Si)\leq \min(\eps_3, \eps_5)/C_2$. 

We can therefore apply corollary \ref{cr-dt} to finish the proof.
\epn

\section{The area is bounded}

In this section, we assume that $(\Si, \III, \Nat)$ is convex (as in
definition \ref{df-cvx}) and has curvature $\Kt\geq K_4$ and torsion
$\|\taut\|\leq \tau_0$, with $4K_4 > \tau_0^2$. We will prove lemma
\ref{borne-aire}, which states that
$(\Si, \III)$ has bounded area. This will be achieved through the
following lemmas:

\blm \label{approx}
For any $\eps>0$, there exists a simply connected domain $\Om\subset
\Si$, $\Omb\subset \Si$, with 
locally convex boundary, such that $d_{\III}(\dr_{\III}\Om,
dr_{\III}\Si)\leq \eps$.
\elm

\blm \label{comp0}
$(\Om, \III)$ can not be complete.
\elm

\blm \label{comp2}
$\dr_{\III} \Om$ can not have a non-compact component.
\elm

\blm \label{comp1}
If $\dr_{\III}\Om$ is a closed curve, then the area of $\Om$ is at most
$2\pi/K_4$. 
\elm

Lemma \ref{comp1} is an immediate consequence of the Gauss-Bonnet
theorem, so that the rest of this section contains the proofs of lemmas
\ref{approx}, \ref{comp0} and \ref{comp2}. Lemma \ref{borne-aire}
follows: by lemma \ref{approx}, any compact subset of $\Si$ should be
contained in a domain $\Om$ with locally convex boundary, which should
have area at most $2\pi/K_4$ by lemmas \ref{comp0}, \ref{comp2} and
\ref{comp1}.

\bpn{of lemma \ref{approx}}
Let $\cE$ be the set of open simply connected domains $\Om\subset \Si$
such that $\dr_{\III}\Om\setminus \dr_{\III}\Si$ is locally convex, and that
$d_{\III}(\dr_{\III}\Om, \dr_{\III}\Si)\leq \eps$. $\cE$ is ordered by
inclusion. Let $\Om_0$ be a minimal element of $\cE$. We want to prove
that $\Omb_0\subset \Si$; we proceed by contradiction, and suppose that
there exists a point $x_0\in \dr_{\III} \Om_0\cap \dr_{\III}\Si$.

Let $x_1\in \Om$ be such that $d(x_0, x_1)\leq \eps_0/2$ in $\Om_0$. Let
$c:[0, L)\rightarrow \Om_0$ be a smooth curve of length $L$
with $c(0)=x_1$ and $\lim_{t\rightarrow L} c(t)=x_0$, with $L\leq t_g$,
where $t_g$ comes from corollaries \ref{cr-ido} and \ref{diff-loc}. 
Let $t_0$ be the
supremum of all $t\in [0, L)$ such that there exists a one-parameter
family $(g_t)_{t\in [0, t]}$ of geodesic segments, with $g_t$ going from
$c(0)$ to $c(t)$. If $t_0<L$, then:
$$ \lim_{t\rightarrow t_0}d_{\III}(g_t, \dr_{\III}\Om)=0~. $$

If $t_0<L$, then, as
$t\rightarrow t_0$, $g_t$ would approach either a point of
$\dr_{\III}\Om_0\setminus \dr_{\III}\Si$ --- but this is impossible
because $\dr_{\III}\Om_0\setminus \dr_{\III}\Si$ is locally convex -- or
a point of $\dr_{\III}\Om_0\cap \dr_{\III}\Si$ -- and this is impossible
by definition \ref{df-cvx}. Therefore, $t_0=L$, and there exists a
geodesic segment $g:=g_L:[0, l)\rightarrow \Om_0$ with $g(0)=x_1$ and
$\lim_{t\rightarrow l}g(t)=x_0$.

For $t\in [0, l)$ and $\theta\in [-\pi, \pi]$, let $\ga_{t, \theta}$ be
the maximal geodesic in $\Si$ with $\ga_{t, \theta}(0)=g(t)$ and
$\angle(g'(t), \ga'_{t, \theta}(0))=\theta$. So $\ga_{t, \theta}$ is a
smooth map from $(-a_{t,\theta}, b_{t,\theta})$ to $\Si$, with
$a_{t,\theta}, b_{t,\theta}\in \R_+^*\cup \{\infty \}$.
Let:
$$ E_{t}:= \{ \theta\in [-\pi, \pi] \; | \; a_{t,\theta}\leq
b_{t,\theta} \}, \; F_{t}:= \{ \theta\in [-\pi, \pi] \; | \; a_{t,\theta}\geq
b_{t,\theta} \}~. $$
For $t>l/2$, $0\in F_t$, while $\pi\in E_t$. Moreover, both $E_t$ and
$F_t$ are closed, so there exists $\theta_t\in E_t\cap F_t$, which means
that $a_{t,\theta_t}=b_{t,\theta_t}\in \R_+\cup\{\infty\}$.

First suppose that there exists a sequence $t_n\rightarrow l$ such that
$a_{t_n, \theta_{t_n}}$ remains bounded.
Then there exists $t$ such that $a_{t, \theta_t}$ is finite and that the
geodesic segments $\ga_{t,\theta_t}$ remains within distance at most 
$\eps/2$ of $\dr_{\III}\Si$ because $(\Si, \III, \Nat)$ is
convex. $\Om_0\setminus \ga_{t,\theta_t}$ has 
at least two connected components, one of which, $\Om_1$, is in
$\cE$: it is convex, and its boundary remains within distance at most
$\eps$ of $\dr_{\III}\Si$. But this contradicts the minimality of
$\Om_0$, and this finishes the proof in this case.

Now suppose that $a_{t,\theta}\rightarrow \infty$. Since $\Si$ has a 
convex boundary, for any $a>0$ and any $\eps'>0$, there exists $t$ close
to $l$ such that $\ga_{t,\theta_t}([-a,a])$ remains within distance
$\eps'$ of $\dr_{\III}\Si$. 
Call $g_0$ the restriction of $\ga_{t,\theta_t}$ to $\R_+$, and apply
proposition \ref{fam-geod}. It shows that, if $a$ is large enough and
$\eps'$ small enough, there exists $\theta>0$ such that $g_{\theta}$
(defined as in proposition \ref{fam-geod})
remains within distance $\eps$ of $\dr_{\III}\Si$ (because it remains
close to $g_0$, point (2) of \ref{fam-geod}) but goes from $\ga_{t,
\theta_t}(0)$ to $\dr_{\III}\Si$ (because of point (3.) of \ref{fam-geod}). 

\vspace{0.3cm}
\centerline{\psfig{figure=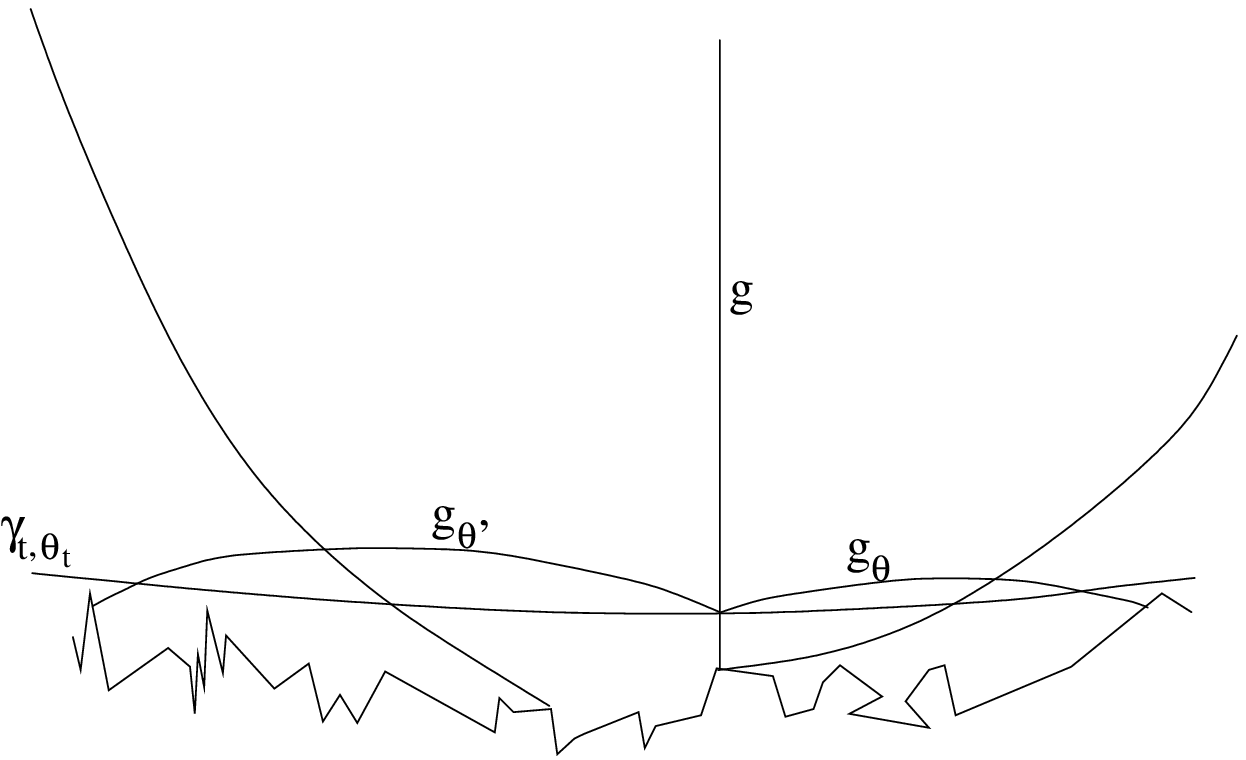,height=5cm}}
\centerline{\bf Figure 7.1} \vspace{0.3cm}

Then do the same for the geodesic ray $g'_0$ defined by $g'_0(s)=\ga_{t,
\theta_t}(-s)$, to obtain another similar geodesic segment
$g'_{\theta'}$ which remains within distance $\eps$ of $\dr_{\III}\Si$,
and goes from $\ga_{t, \theta_t}(0)$ to $\dr_{\III}\Si$. 

The rest of the proof can be done as in the case where $a_{t, \theta_t}$
remains bounded, because $g_{\theta}\cup g'_{\theta'}$ ``cuts'' a convex
domain in $\Om$ whose boundary remains within distance $\eps$ of
$\dr_{\III}\Si$, a contradiction.
\epn

\bigskip

We now come back to deformations of convex curves, as in the previous
section. The following is a consequence of proposition \ref{edo}.

\bprop \label{edo7}
For any $\eps>0$ small enough, there exists $M'_1$ as follows. Choose
$N_1>0$. For any
$u\in C^\infty(\R, [-1/\eps, 1/\eps])$, there exists $y\in
C^0([0, s_0], [0, M'_1])$ which:
\begin{enumerate}
\item is piecewise $C^\infty$; 
\item vanishes outside
$[-N_1-M'_1, N_1+M'_1]$;
\item is at least $1$ in $[-N_1, N_1]$; 
\item satisfies (as a distribution): 
\beq \label{ido1} 
y'' \geq (yu)'-(\eps+\frac{u^{2}}{4}) y~;
\eeq
\item is $M'_1$-Lipschitz.
\end{enumerate}
\eprop

\bpv
Let $x_0:=-N_1$. 
First apply proposition \ref{edo} after translating the
origin to $-N_1$ in $\R$, so as to obtain $x_1\in [-N_1+S_0, -N_1+S_1]$
and a solution $y_1:[x_0, x_1]\rightarrow [1, M_0]$ of (\ref{ido1})
with:
$$ y_1(x_0)=y_1(x_1)=1 ~~ \mbox{and} ~~ y_1'(x_0)=u(x_0)+4, ~
y_1'(x_1)\leq u(x_1)+4~. $$   
Apply proposition \ref{edo} again, now after a translation of the origin
to $x_1$; this provides us with $x_2\in [x_1+S_0, x_1+S_1]$ and with a
solution $y_2:[x_1,x_2]\rightarrow [1, M_0]$ of (\ref{ido1}) with:
$$ y_2(x_1)=y_2(x_2)=1  ~~ \mbox{and} ~~ y_2'(x_1)=u(x_1)+4, ~ y_2'(x_2)\leq
u(x_2)+4~. $$
Repeat this procedure to find a sequence $x_3 \leq \cdots \leq x_N$ with
$x_{k+1}-x_k\in [S_0, S_1]$ and $x_{N-1}\leq N_1<x_N$, and functions
$y_k:[x_{k-1}, x_k]\rightarrow [1, M_0]$ which are solutions of
(\ref{ido1}) with:
$$ y_k(x_{k-1})=y_k(x_k)=1  ~~ \mbox{and} ~~ y_k'(x_{k-1})=u(x_{k-1})+4,
~ y_k'(x_k)\leq u(x_k)+4~. $$

Apply proposition \ref{edo} once more to find $x_{N+1}\in [x_N+S_0,
x_N+S_1]$ and a solution $y_{N+1}:[x_N, x_{N+1}]\rightarrow [0, M_0]$ of
(\ref{ido1}) with: 
$$ y_{N+1}(x_{N})=1~, ~~ y_{N+1}(x_{N+1})=0~, ~~
y_{N+1}'(x_{N})=u(x_{N})+4~, ~~ y_{N+1}'(x_{N+1})\leq 0~. $$
Do the same to find $x_{-1}\in [x_0-S_1,
x_0-S_0]$ and a solution $y_{-1}:[x_N, x_{N+1}]\rightarrow [0, M_0]$ of
(\ref{ido1}) with: 
$$ y_{-1}(x_{0})=1~, ~~ y_{-1}(x_{-1})=0~, ~~
y_{-1}'(x_0)=0~, ~~ y_{-1}'(x_{-1})\geq 0~. $$

Now define $y:\R\rightarrow \R_+$ as the function whose restriction to
$[x_k, x_{k+1}]$ is $y_{k+1}$ for $-1\leq k\leq N$, and which is zero
outside $[x_{-1}, x_{N+1}]$. Note that, if $M'_1:=2S_1$, then $y$
vanishes outside $[-N_1-M'_1, N_1+M'_1]$. Moreover, it is clear that $y$
is a (weak) solution of (\ref{ido1}).
\epv

The previous proposition provides the tool needed to deform convex
curves while increasing their curvature.

\bcr \label{dfo-R}
There exist $c_1>0$ as
follows. Let $\kappa_m\in \R_+$, and let $g_0:\R\rightarrow
\Si$ be a smooth, convex curve parametrized at speed one, with
curvature $\kappa\geq \kappa_m$ as a measure. Let $N_1>0$. Then there 
exists $T>0$ and a deformation $(g_t)_{t\in [0, T]}$ such that, for
each $t\in [0, T]$:
\begin{enumerate}
\item $g_t$ is a convex curve with curvature $\kappa\geq \kappa_m$, and
$\kappa\geq \kappa_m+t c_1$ along $g_t([-N_1, N_1])$; 
\item $g_t\equiv g_0$ outside $[-N_1-M'_1, N_1+M'_1]$;
\item for each $s\in \R$, either $(\dr_t g_t)(s)$ is zero, or its
orthogonal is a support direction of $g_t$, and its norm is at most
$M'_1$. 
\end{enumerate}
\ecr

\bpv
Let $(g_n)_{n\in \N^*}$ be a sequence of smooth curves, $g_n:\R\rightarrow
\Si$, such that:
\begin{itemize}
\item $\forall s\in \R, \lim_{n\rightarrow \infty} g_n(s)=g_0(s)$;
\item for $n\leq m$, $g_n$ lies entirely on the concave side of $g_m$;
\item the curvature of $g_n$ is at least $-\alpha_n<0$, where
$\lim_{n\rightarrow \infty}\alpha_n=0$.
\end{itemize}
The existence of such an approximating sequence is not too difficult to
prove. The $(g_n)$ are not parametrized at speed one.

We suppose (without loss of generality) that, for $n\in \N^*$ and $s\in
\R$, $Jg'_n(s)$ is towards $g_0$. For $n\in \N^*$ and
$s\in \R$, let:
$$ u_n(s) := \tau(J g'_n(s))~. $$
Apply corollary \ref{edo7} to obtain a sequence of piecewise smooth,
$M'_1$-Lipschitz functions $(y_n)_{n\in \N^*}$ with:
$$ y_n''\geq (y_n u_n)' - \left(\eps+\frac{u_n^2}{4}\right) y_n~, $$
with $y_n(s)=0$ when $s\not\in [-N_1-M_1', N_1+M_1']$ and $y_n(s)\in [1,
M'_1]$ when $s\in [-N_1, N_1]$. Since the $y_n$ are Lipschitz, we can (by
taking a subsequence) suppose that they are $C^0$-converging to a
Lipschitz function $y:\R\rightarrow [0, M'_1]$.

Let $T'\in \R_+\cup\{ \infty\}$ be the largest $t$ such that, for each
$n\in \N^*$ and each $s\in [-N_1, N_1]$, $\exp^{\Nat}_{g_n(s)}(tJg'_n(s))$
is defined. Then $T'>0$ by compactness. For $n\in \N^*$ and $s\in [0,
T')$, let:
$$ h_{n,t}(s) := \exp^{\Nat}_{g_n(s)}(tJg'_n(s))~. $$
For $n\in \N^*$ and $t\in [0, T')$, $h_{n,t}$ is a curve which might not
be embedded, but which, for $t$ small enough, is immersed. It
differs from $g_n$ only in $[-N_1-M_1', N_1+M_1']$. Moreover, corollary
\ref{cr-df} and a simple compactness argument show that there exist
$N\in \N^*$, $T\in (0, T')$ and $c>0$ such that the curvatures
$\kappa_{n, t}$ of the curves $h_{n, t}$ satisfy:
\beq \label{ineg-hnt}
\forall n\geq N, \forall t\in [0, T), \forall s\in \R,
\kappa_{n,t}(s)\geq \kappa_m + c t -\eps_n~,
\eeq
where the left-hand side is a measure and $\eps_n\rightarrow 0$.

\vspace{0.3cm}
\centerline{\psfig{figure=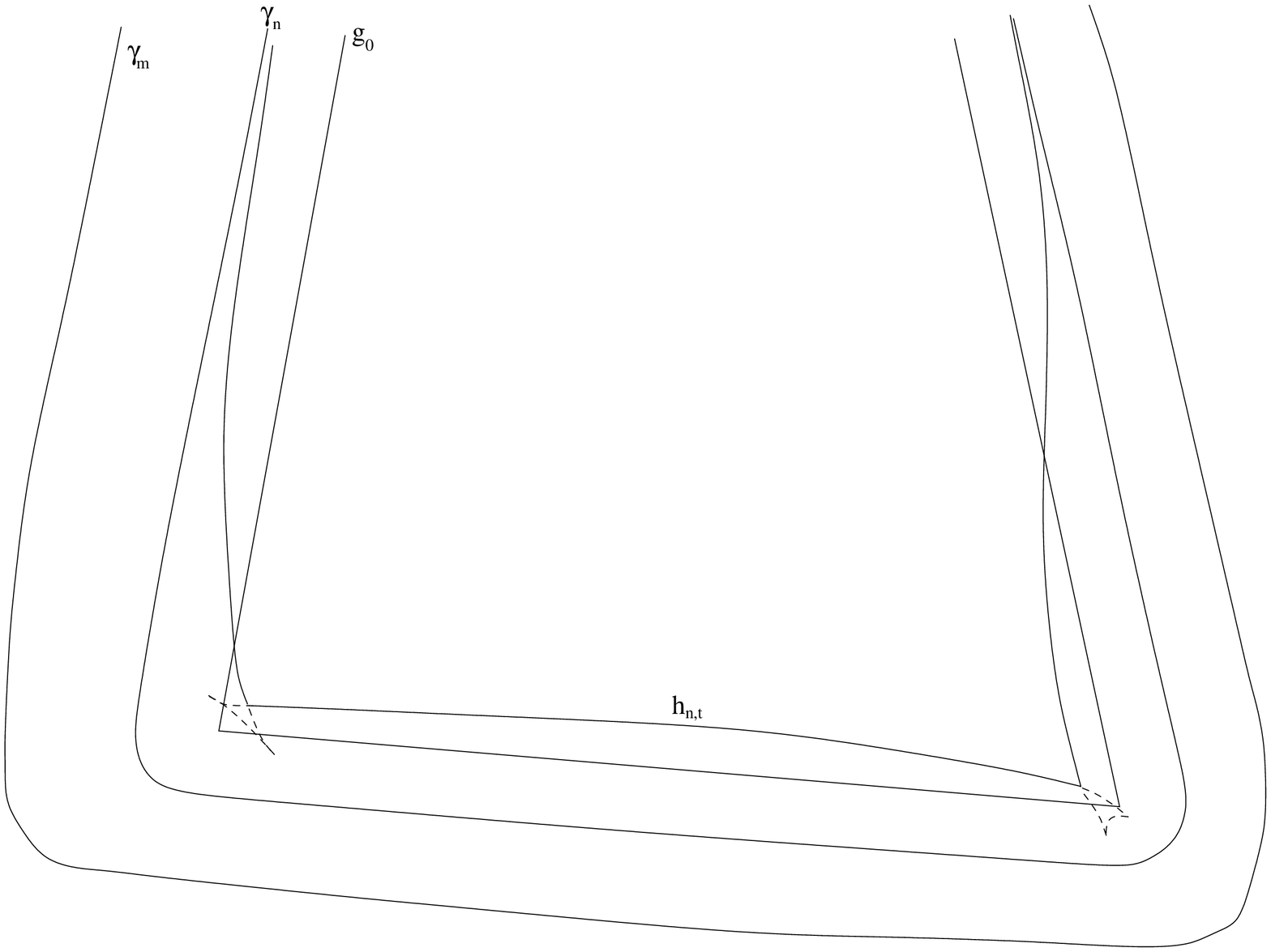,height=7cm}}
\centerline{\bf Figure 7.2} \vspace{0.3cm}

Since the $h_{n,t}$ are curves and differ from $g_n$ only in a
compact set, they separate $\Si$ into several connected components, two
of which are non compact; we call $\Om_{n,t}$ the non-compact connected
component of $\Si\setminus h_{m,t}$ whose intersection with the concave
side of $g_n$ is (empty or) compact. Equation (\ref{ineg-hnt}) shows
that the boundary of 
$\Om_{n,t}$ in $\Si$ is locally convex, with curvature at least
$\kappa_m+c t-\eps_n$. 
 
Finally, for $t\in [0, T)$, set:
$$ g_t := \dr\left(\cap_{n=0}^\infty \Om_{n,t}\right)~. $$
It is not difficult to check that $(g_t)$, with an adequate
parametrization, satisfies the conclusion of corollary \ref{dfo-R}.
\epv

As a consequence, the same kind of deformation can be done not only for
small $t$, but for all $t$:

\bcr \label{dfo-Rt}
Let $\Om$ be a closed subset of $\Si$ with locally convex
boundary. Suppose that some connected component of $\dr\Om$ is a
complete, non compact curve, parametrized at speed one by
$c_0:\R\rightarrow \dr\Om$. Choose $N_1>0$. Then there exists a
deformation $(c_t)_{t\in \R_+}$ such that, for
each $t\in \R_+$:
\begin{enumerate}
\item $c_t$ is a convex curve in $\Om$, with curvature $\kappa\geq
\kappa_m$, and 
$\kappa\geq \kappa_m+t c_1$ along $g_t([-N_1, N_1])$; 
\item $g_t\equiv g_0$ outside $[-N_1-M'_1, N_1+M'_1]$;
\item for each $s\in \R$, either $(\dr_t g_t)(s)$ is zero, or its
orthogonal is a support direction of $g_t$, and its norm is at most
$M'_1$. 
\end{enumerate}
\ecr

\bpv
The underlying idea is to apply corollary \ref{dfo-R} recursively, to
obtain the existence of 
such a deformation for $t\in [0, T]$ for some $T>0$. The formal proof,
however, has to be done in a slightly different way. Suppose that such a
deformation can not exist for all $t\in \R_+$. Let $E$ be the set of
couples $(t, (g_s)_{s\in [0, t)})$, where $t>0$ and $(g_s)_{s\in [0,t)}$
satisfies the conditions demanded, but only until time $t$. 

There is a natural order on $E$, with:
$$ (t, (g_s)_{s\in [0, t)}) \leq (t', (g'_s)_{s\in [0, t')}) $$
if $t\leq t'$ and $g_s=g'_s$ for $s\leq t$. $E$ has a maximal element,
say $(t_0, (g^0_s)_{s\in [0, t_0)})$. For $u\in \R$, let:
$$ g^0_{t_0}(u) := \lim_{t\rightarrow t_0} g^0_t(u) \in \Sib~. $$
Because of the convexity of the $(g_s^0)$ and because
$\dr_{\III}\Si$ has no concave point, $g_{t_0}$ is a convex curve. Thus
one can apply corollary \ref{dfo-R} to $g_{t_0}$, and this contradicts
the maximality of $(t_0, (g^0_s)_{s\in [0, t_0)})$.
\epv

We now consider related questions for deformations of curves which are
topologically $S^1$, and which are not necessarily convex, but have
curvature bounded from below.

\bcr \label{dfo-S}
There exist $c_2>0$ and $M'_2\geq 1$ as
follows. Let $\kappa_m\in \R$, $L\in \R_+^*$, and let
$g_0:\R/L\Z\rightarrow \Si$ be a curve parametrized at speed one,
with curvature $\kappa\geq \kappa_m$ (with the normal oriented towards
the non-compact side of $g_0$). Then there 
exist $T_2>0$ and a deformation $(g_t)_{t\in [0, T_2]}$ such that, for
each $t\in [0, T_2]$:
\begin{enumerate}
\item $g_t$ is a curve which bounds a compact set, with curvature
$\kappa\geq \kappa_m+t c_2$;  
\item for each $s\in [0, L]$, $\| (\dr_t g_t)(s)\| \in [1, M'_2]$. 
\end{enumerate}
\ecr

\bpv
It is similar to the proof of corollary \ref{dfo-R}. First choose a
sequence of smooth curves $g_n:\R/L\Z$ converging to $g_0$, such that
$g_n$ is in the interior of $g_m$ for $n\leq m$ and that the curvature
$\kappa_n$ of $g_n$ is at least $\kappa_m-\alpha_n$, with
$\lim_{n\rightarrow\infty} \alpha_n=0$.

For $n\in \N^*$ and $s\in \R/L\Z$, let:
$$ u_n(s) := \tau(Jg'_n(s))~, $$
where we suppose again that $Jg'_n$ is towards the non-compact side of
$g_n$. Let $\ut_n$ be the lift of $u_n$ to a function on $\R$. Apply
proposition \ref{edo7} to $\ut_n$, to obtain a piecewise smooth,
Lipschitz function $\yt_n:\R\rightarrow M'_1$ which vanishes outside
$[-L-M'_1, L+M'_1]$ and which is at least $1$ on $[-L,L]$. Let
$y_n:\R/L\Z\rightarrow \R$ be the function defined by:
$$ y_n(u) := \sum_{s\in u} \yt(s)~, $$
where only finitely many terms are non-zero by definition of
$\yt_n$. After multiplying it by a constant, $y_n$ is a Lipschitz,
piecewise smooth 
function from $\R/L\Z$ to $[1, M'_0]$ (for some $M'_0>1$ which depends
on $L$ and on $M'_1$) which is a solution of (\ref{edo1}). 

The rest of the proof can be done quite like in the proof of corollary
\ref{dfo-R}, so we leave the details to the reader.
\epv

\bcr \label{dfo-SR}
Corollary \ref{dfo-S} is true for any $T_2>0$.
\ecr

\bpv
Like the proof of corollary \ref{dfo-Rt} from corollary \ref{dfo-R}.
\epv

We now have enough results on the deformations of curves, and we turn to
another simple property: a convex, complete curve which separates a
convex subset of $\Si$ into two parts can not be ``too'' curved.

\bprop \label{prop-fin}
There exists a constant $\kappa_{0}(K_{4},\tau_{0})$ as follows. Let
$\Om$ be a closed, locally convex subset of $\Si$, and let
$\rho:\R \rightarrow \Om$ be a convex, injective curve, parametrized at
speed one, which separates $\Om$ into two connected components, and with
curvature: 
\beq \label{dec-courb}
\kappa = \kappa_{1}+\kappa_{m}~,
\eeq
where $\kappa_{1}>0$ is a constant and $\kappa_{m}$ is a positive
measure. Then $\kappa_{1}\leq \kappa_{0}$.
\eprop

\bpv
First note that, by a direct approximation argument, it is enough to
prove the result when $\rho$ is smooth, so we suppose that is the
case. Let $t_0\in \R$.
By corollary \ref{ex-geod}, there exists $\eps>0$ (depending only on
$K_{5}$ and $\tau_{0}$), 
such that $\exp_{\rho(t_0)}^{\Nat}$ is a
diffeomorphism from the subset of the ball of radius $\eps$ where it is
defined onto its image. 
Therefore, for all $t\in [t_{0}, t_{0}+\eps]$, there exists a unique
$\Nat$-geodesic  $\ga_{t}:[0,L_{t}]\rightarrow \Om$
of minimal length between $\rho(t_{0})$ et $\rho(t)$. 
For $t\in [t_{0}, t_{0}+\eps]$, let $\theta_{1}(t)$
be the angle between $\rho'(t_{0})$ and $\ga_{t}'(0)$, $\theta_{2}(t)$
the angle between $\ga'_{t}(L_{t})$ and $\rho'(t)$, $D(t)$ the domain in
$\Om$ bounded by $\rho([t_{0}, t])$ and by $\ga_{t}$, and $A(t)$ its area.

From the Gauss-Bonnet theorem, for each $t\in [t_{0}, t_{0}+\eps]$:
\begin{eqnarray}
\theta_{1}+\theta_{2} & = & \int_{t_{0}}^{t_{1}}\kappa(ds) + \int_{D(t)}
\Kt da \nonumber \\
& \geq & \kappa_{1}(t-t_{0}) + K_{4} A(t)~. \nonumber
\end{eqnarray}
But it is easy to check, using equation (\ref{ido-y}) of corollary
\ref{cr-ido}, that, if $L_{t}$
remains small enough (so that $\eps$ remains smaller than a constant
depending only on $K_{5}$ and on $\tau_{0}$), then $A(t)$ is bounded by:
\begin{equation} \label{maj-At}
A(t) \leq \int_{t_{0}}^{t} L(s)ds \leq 2\int_{t_{0}}^{t} (s-t_{0}) ds
\leq (t-t_{0})^{2}~.
\end{equation}
As a consequence of those two equations, if $\kappa_{1}$ is large enough
(larger than a constant depending only on $K_{5}$ and on $\tau_{0}$),
there exists $t_{1}\in
[t_{0}, t_{0}+\eps]$ such that:
$$
\theta_{1}(t_{1})+\theta_{2}(t_1)=2\pi~.
$$
Then there exists $t_{2}\in [t_{0},t_{1}]$ verifying one of the
following properties:
\begin{enumerate}
\item either $\theta_{1}(t_{2})=\pi$;
\item or $\theta_{1}(t_{2})\leq \pi$ and $\theta_{2}(t)=\pi$.
\end{enumerate}

\vspace{0.3cm}
\centerline{\psfig{figure=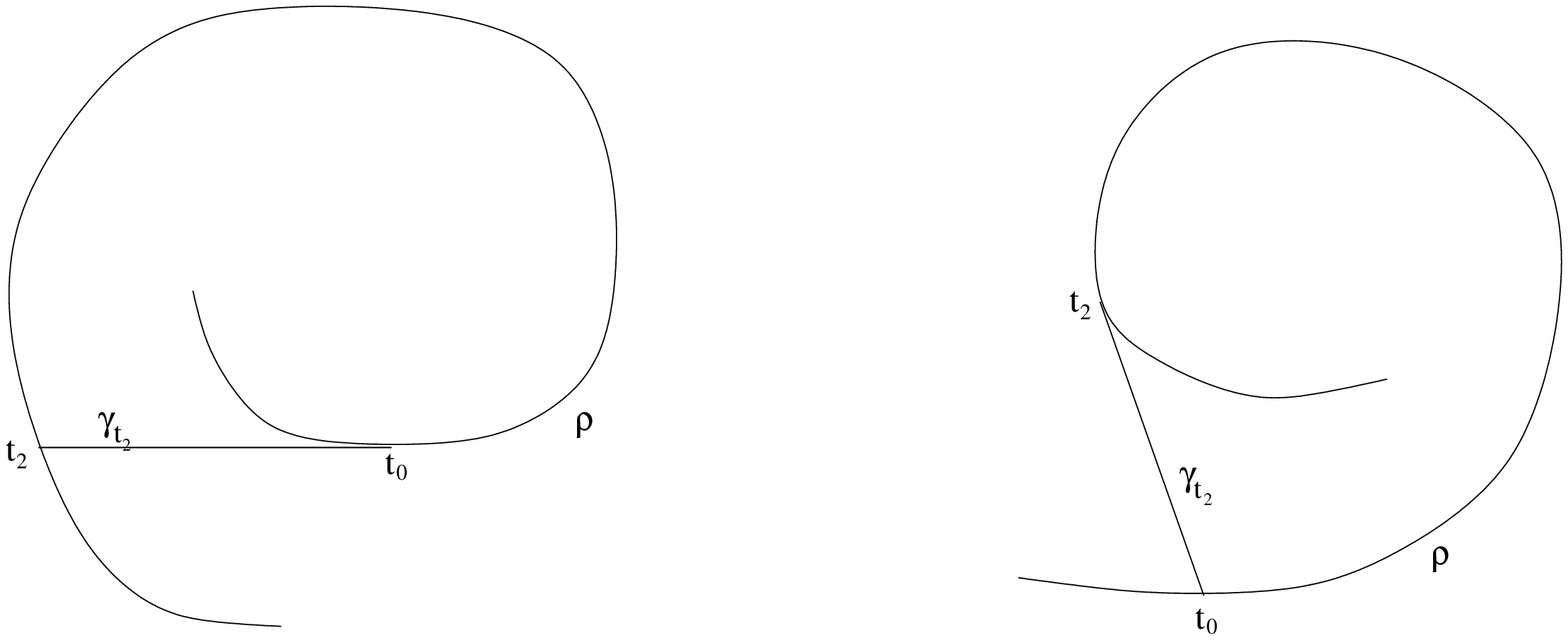,height=4cm}}
\centerline{\bf Figure 7.3} \vspace{0.3cm}

Recall that $\rho$ is injective; therefore, using
(\ref{dec-courb}) and the upper bound (\ref{maj-At}) on $A(t_{2})$, one
sees that, if $\eps$ is small enough (smaller than a constant which this
time depends on $K_{4}$ and on $\tau_{0}$), then:
\begin{itemize}
\item in the first case, that $\rho(]-\infty, t_{0}])$ remains in the
domain of $\Om$ bounded by $\rho([t_{0}, t_{2}])$ and by $\ga_{t_2}$;
\item in the second case, that $\rho([t_{2}, \infty])$ remains in the
domain of $\Om$ bounded by $\ga_{t_2}$.  
\end{itemize}

In both cases, ``half'' of $\rho$ remain in a compact domain of $\Om$,
and therefore $\rho$ can not separate $\Om$ in two parts.
\epv

We can show that $\Om$ can not be complete:

\bigskip

\bpn{of lemma \ref{comp0}}
Choose a smooth, simple closed curve $\ga_{0}$ in $\Om$, 
with its unit normal oriented towards the non-compact connected
component of $\Om\setminus \ga_{0}$ (there exists one because $\Om$ is
complete and simply connected). Then the geodesic curvature of $\ga_{0}$
can be written as:
\beq \label{kappa}
\kappa = -\kappa_{1}+\kappa_{m}~,
\eeq
with $\kappa_{1}\in \R_{+}$ and $\kappa_{m}\geq 0$. 
Apply corollary \ref{dfo-SR}, to obtain a continuous family
$(\ga_{t})_{t\in \R_{+}}$ of curves with curvature bounded
below by $-\kappa_1+c_2 t$. For $t$ large enough, this contradicts
proposition \ref{prop-fin}.
\epn

\bigskip

Finally, we can now prove that $\dr \Om$ can not contain any non-compact
curve.

\bigskip

\bpn{of lemma \ref{comp2}}
Suppose that $\dr_{\III}\Om$ contains a non compact curve $g_0$. Since
$\Om$ is simply connected, proposition \ref{prop-gagat} shows that $g_0$ is
at non-zero distance from the other connected components of
$\dr_{\III}\Om$. 

Corollary \ref{dfo-Rt} thus shows that there exists a continuous
deformation $(g_t)_{t\in [0, T)}$ of $g_0$ in $\Om$, which goes on until
time $T$ with either $T=\infty$ or $\lim_{t\rightarrow T}d_{\III}(g_t,
\dr\Om\setminus g_0)=0$. But proposition \ref{prop-gagat} excludes this
possibility, so that $T=\infty$. This contradicts
proposition \ref{prop-fin} just as in the proof of lemma \ref{comp0}
above. 
\epn

\bigskip

It is a natural question whether the hypothesis we had to make
concerning the relationship between $\tau_0$ and $K_4$ are really
necessary to obtain a bound on the area of $(\Om, \III)$. This is all the
more important since the main geometric result of this paper, theorem
\ref{thm-a}, is limited precisely by this hypothesis.

The following example shows that this relation is in a sense
optimal. Note that, for homogeneity reasons, any relationship between
the curvature and the torsion should relate the curvature to the square
of the torsion.

Consider the hyperbolic plane $H^{2}$, with the connection 
$\Na_{t}$ obtained by adding to the Levi-Civita connection
$\Na_{0}$ a 1-form $\beta_{t}$, so that:
$$
\Na^{t}_{X}Y = \Na^{0}_{X}Y +\beta_{t}(X)JY~,
$$
with $\beta_{t}=-t u_{\theta}^{*}$, where $u_{\theta}$ is at each point
the unit normal bundle to the geodesic coming from $0$ (with the usual
orientation) and $u_{\theta}^{*}$ is the dual 1-form.

It is easy to check the following points:
\begin{itemize}
\item the torsion $\tau_{t}$ of $\Na_{t}$ is such that: $\| \tau_{t}
\| = t$;
\item for $\beta>1$ and $r$ large enough, the complement
$\Ga_{r}$ of the ball centered at $0$ of radius $r$ is
$\Na_{t}$-convex, with infinite area;
\item the curvature of $\Na_{t}$ is: $K_{t}=t.\mbox{th}(r)-1$ at
distance $r$ from $0$.
\end{itemize}

The lower bound on $K_{t}$ is therefore $K^{m}_{t}=t-1$,
and the smallest value of $\tau_{t}^{2}/K^{m}_{t}$ is obtained for
$t=2$, it is equal to $\frac{1}{4}$. So, using the same notations as
above, for $K_{4}=\tau_{0}^{2}/4$, there is already a counter-example to
the results proved above for $K_{4}>\tau_{0}^{2}/4$.

An interested reader might check that one can built other examples,
based on deformations of the Levi-Civita connection of a positively
curved surface (e.g. an annulus in $S^2$ with its canonical metric) such
that convex 
domains with infinite area, for instance ``strips'' between two convex
curves with constant curvature, exist. But the limiting relations between
$K_4$ and $\tau_0$ are the same as in the hyperbolic-based example
above.

\section{Some examples and further statements}

This section contains some other, more precise results like theorem
\ref{thm-a}, and also some examples which indicate that theorem
\ref{thm-a} is, in some sense, optimal. 

First note that the proof which was given actually shows a little more
than what was stated, namely:

\btm \label{th1}
Let $K_1, K_2, K_3\in \R$ be such that $K_1<0$ and that $K_1<K_2\leq K_3$.
Let $(\Si, \si)$ be a complete Riemannian surface, and let $(M,\mu)$ be
a Riemannian 3-manifold. Suppose that the curvature $K_\Si$ of $\Si$ is
bounded above by $K_1$, and that, for all $m\in M$, the maximal and
minimal curvatures of the 2-planes in $T_mM$, $K_M$ and $K_m$,
are in $[K_2,K_3]$ and such that:
\beq \frac{K_M-K_m}{2\sqrt{(K_1-K_M)(K_1-K_m)}} \leq \tau_0~, 
\label{cond0} \eeq
with
\beq
 \tau_0^2 < \frac{4K_1}{K_1-K_3} ~ \mbox{if} ~ K_3\leq 0 ~ , ~
\tau_0^2 < 4 ~ \mbox{if} ~ K_3\geq 0~.
\label{tau0} \eeq
Suppose further that the gradient of the sectional curvature of
$(M,\mu)$ (on $M$) and the gradient of $(K_{\si})^{-1/2}$ (on $\Si$) are
bounded. Then there exists no $C^{3}$ isometric immersion of $(\Si,\si)$
into $(M,\mu)$.
\etm

Note that the precise value of $K_2$ plays no role in this statement;
but we need to know that $K_2>K_1$ to obtain an upper bound on the
curvature of $\Nat$. This hypothesis is related to the  ``uniform
hyperbolicity'' of the immersion.

\bigskip

Theorem \ref{thm-a} is strongly related to a well-knowsn family of PDEs,
the Monge-Amp{\`e}re PDEs of hyperbolic type. Those are usually written, on
a domain $\Om \subset \R^{2}$, as:
\beq \label{mah}
\frac{\dr^{2} u}{\dr x^{2}} \frac{\dr^{2} u}{\dr y^{2}} - \frac{\dr^{2}
u}{\dr x \dr y} = -b~,
\eeq
where $b$ might be a function on $\Om$, or depend on $u$ and maybe of
its first derivatives. 

Equation (\ref{mah}) can be written as an equation on the bundle
morphism $H$ associated to the hessian
$\Na^{2}u$ of $u$ (by: $ (\Na^{2}u)(X,Y) =
\langle HX | Y \rangle = \langle X | HY \rangle$)
with first order conditions meaning that $H$ is the hessian of a
function; we get thus:
\beq \label{mah-g}
\left\{ \begin{array}{c}
\det(H)=-b \\
d^{\Na}H = 0~,
\end{array} \right. 
\eeq
This form seems well adapted to the study of Monge-Amp{\`e}re equations over
Riemannian surfaces other than $\R^2$. 
in particular, if $b$ is everywhere equal to $K_{\si}-K_{0}$, then this
equation is verified by the second fundamental forms of the immersions
into space-forms with curvature $K_{0}$. For immersions into a
3-manifold with non-constant curvature, a right-hand side appears in the
first equation of (\ref{mah-g}). This generalization is not equivalent
to what is obtained by writing (\ref{mah}) in term of the hessian of
$u$; if $(M,\mu)$ is a Riemannian manifold with Levi-Civita connection
$\Na$ and curvature tensor $R$, if $u:M\rightarrow \R$ is a $C^{2}$
function, and if $H$ is the bundle morphism associated to its hessian, then:
$$
\forall m\in M, \; \forall X,Y\in T_{m}M, \; (d^{\Na}H)(X,Y) = -
R(X,Y)(Du)~, 
$$ 
which differs from (\ref{mah-g}) because the 1-jet of $u$ appears.

One can check that the proof given for theorem \ref{thm-a} also proves
the following:

\btm \label{th-mah}
Let $\eps_{0}>0$, $0<b_{m}\leq b_{M}$ and $\tau_{0}\geq 0$ be such that 
$b_{M}\tau_{0}^{2} < 4\eps_{0} b_{m}^{2}$; let $(\Si,\si)$ be a complete
Riemannian surface with curvature $K\leq -\eps_{0}$,
$b:\Si\rightarrow [b_{m},b_{M}]$ be a $C^{1}$ function with 
bounded gradient on $\Si$, and let $\tau$ be a $C^{0}$ vector field on
$\Si$ with $\|\tau\|\leq \tau_{0}$. Then the system:
\beq
\left\{ \begin{array}{c}
\det(H)=-b \\
d^{\Na}H = \tau \otimes \nu_{\si}
\end{array} \right. 
\eeq
(where $\nu_{\si}$ is the area form associated to $\si$ and $H$
is a symmetric endomorphism field on $\Si$) has no 
$C^{1}$ solution on $\Si$.
\etm

The point is that such a solution $H$ would allow the definition of a
``virtual third fundamental form'' on $\Si$ as:
$\III(X,Y)=\si(HX,HY)$. In addition, we could define a connection $\Nat$
(as in section 2) compatible with $\III$, with curvature 
$\Kt=-K/b$, and torsion $\taut=-H^{-1}\tau/b$. With the hypotheses of
theorem \ref{th-mah}, we would have:
$$ \|\taut\|_{\III} = \|\tau\|_{\si}/b \leq \tau_{0}/b_{m}~, $$
so that:
$$
\|\taut\|_{\III}^{2} \leq \frac{\tau_{0}^{2}}{b_{m}^{2}} <
\frac{4\eps_{0}}{b_{M}} \leq 4\Kt~,
$$ 
and the analog of lemma \ref{lem-e} would indicate that
$(\Si,\III,\Nat)$ should be convex; while the analog of lemma
\ref{borne-aire} would lead to a contradiction. 

\bigskip

There are also various possible improvements of theorem \ref{thm-a}. For
instance, almost all the proof takes place ``at infinity'', while the
interior of $(\Si, \III, \Nat)$ is important essentially only in the
proof of lemma \ref{borne-aire}. Therefore, one can check that, if
$(\Si, \si)$ is simply connected but satisfies the hypothesis of theorem
\ref{thm-a} only outside a compact set, then isometric immersions remain
impossible.

\bigskip

It is natural to wonder to what extend the conditions in theorem 
\ref{thm-a} are really necessary. It is not clear concerning the
hypothesis that the gradient of the sectional curvatures of $(\Si,\si)$
and of 
$(M,\mu)$ are bounded, but the following example shows that the
inequalities in theorem \ref{thm-a} are necessary.

Let $g_{\lambda}$ be the symmetric 2-form  defined on $\R^{3}$ as:
$$
g_{\lambda} = (1+2\lambda z) \cosh^{2}(y) \cosh^{2}(z) dx\otimes dx +
(1-2\lambda z) 
\cosh^{2}(z) dy \otimes dy + dz\otimes dz~,
$$
for $\lambda>0$. It is a Riemannian metric in a neighborhood of
$P_{0}=\{(x,y,z)\in \R^{3} \; | \; z=0 \}$. Let 
$V_{\eps}$ be such a neighborhood:
$$
V_{\eps}=\{(x,y,z)\in \R^{3} \; | \; |z|\leq \eps \}~.
$$
When $\lambda=0$, $g_{\lambda}$ is hyperbolic (i.e. it has constant
curvature $-1$.

$g_{\lambda |P_{0}}$ is hyperbolic, so we have an isometric embedding of
$H^{2}$ into a Riemannian manifold with boundary  $(V_{\eps}, g_{\lambda
|V_{\eps}})$. 
 
A rather boring computation leads to the following expression of the
Riemann curvature tensor of $g_{\lambda}$: if 
$(e_{1},e_{2},e_{3})$ is a orthonormal moving frame made of vectors
directed by $\dr /\dr x$, $\dr /\dr y$ and $\dr /\dr z$, then
\begin{eqnarray} \nonumber
g_{\lambda}(R_{\lambda}(e_{1},e_{2})e_{1},e_{2}) & = & \lambda^{2}-1 \\ 
\nonumber g_{\lambda}(R_{\lambda}(e_{1},e_{3})e_{1},e_{3}) & = &
\lambda^{2}-1 \\ 
\nonumber g_{\lambda}(R_{\lambda}(e_{3},e_{2})e_{3},e_{2}) & = &
\lambda^{2}-1 \\ 
\nonumber g_{\lambda}(R_{\lambda}(e_{1},e_{2})e_{1},e_{3}) & = &
2\lambda \tanh(y) \\ 
\nonumber g_{\lambda}(R_{\lambda}(e_{2},e_{1})e_{2},e_{3}) & = & 0 \\
\nonumber g_{\lambda}(R_{\lambda}(e_{3},e_{1})e_{3},e_{2}) & = & 0~. 
\end{eqnarray}  
It is then a matter of computation to check that, at the point
$(x,y,z)$, the eigenvalues of $R_{\lambda}$ remain in the interval
$[\lambda^{2}-1-2\lambda \tanh(y), \lambda^{2}-1+2\lambda \tanh(y)]$. 
Therefore, for the isometric embedding of $H^{2}$ into $(V_{\eps},
g_{\lambda |V_{\eps}})$ which we have obtained:
$$ 
K_{1}=-1, \; K_{2}=\lambda^{2}-1-2\lambda, \;
K_{3}=\lambda^{2}-1+2\lambda~, 
$$
so that, for $\lambda$ large enough:
$$
(K_{3}-K_{2})^{2}=16\lambda^{2}
$$
and
$$
16(K_{2}-K_{1})|K_{1}|=16\lambda^{2}-32\lambda~.
$$
As $\lambda\rightarrow \infty$, we get very chose to the inequalities in
theorem \ref{thm-a}.
This example is actually related to the one built in section 7.

\bigskip

The condition that 
$K_{1}<0$ is also necessary, because of a classical example:
$T^{2}=\R^2/\Z^2$ admits an isometric embedding into $S^{3}$ with its
canonical metric (up to a factor $\sqrt{2}$. This is because 
$T^{2}$ can be obtained as:
\beq \label{df-tore}
T^{2}=\{(x,y,z,t)\in \R^{4} \; | \; x^{2}+y^{2}=z^{2}+t^{2}=\frac{1}{2}\}
\eeq
with the induced metric, while $S^{3}$ is:
$$
S^{3}=\{(x,y,z,t)\in \R^{4} \; | \; x^{2}+y^{2}+z^{2}+t^{2}=1\}
$$ 
and the embedding of $T^2$ in $S^3$ follows.

\bigskip

If the target manifold is Lorentzian instead of Riemannian, the
proof given in this paper also applies, with the necessary changes in
the hypothesis concerning the curvature of $\Si$ and of $M$. In fact,
the hypothesis which are needed in this case are such that the only
possible target manifolds are those with constant negative
curvature. This leads to the following result, which was already given
in \cite{efj}:

\btm \label{thL}
Let $\eps\in ]0,-1/2[$, and let $(\Si,\si)$ be a complete Riemannian
surface with curvature $K$ between
$-1+\eps$ and $-\eps$, and such that the gradient $K$ is bounded. Then
there exists no $C^{3}$ isometric immersion of $(\Si,\si)$ into the
anti-de Sitter space $H^3_1$.
\etm

On the other hand, the following assertion is easy:

\bprop
Let $(\Si,\si)$ be a complete Riemannian surface whose curvature $K$ is
larger than some $\eps>0$. There is no isometric immersion $\phi$ of
$(\Si,\si)$ into a Lorentzian manifold $(M,\mu)$ such that, at each
$s\in \Si$:
$$
K^{\Si}(s) > K^{M}(\phi_{*}(T_{s}\Si))~.
$$ 
\eprop  

As a consequence, there is no strictly hyperbolic isometric immersion of
a complete surface into the Minkowski or the de Sitter space.

\section*{Acknoledgements}

This works owes much to many comments and remarks from Fran{\c c}ois
Labourie, and to his encouragements.

Part of this paper was written at the F.I.M. of the E.T.H., Z{\"u}rich; I
would like to express my gratitude for the excellent working conditions
enjoyed there.

\end{document}